%% file: __surfaces_and_area.tex
\documentclass[12pt,oneside]{book}
\usepackage{geometry}                		
\usepackage{graphicx}									
\usepackage{amssymb}\vspace{0.25cm}
\usepackage{mathtools}
\usepackage{setspace}
\usepackage{tikz-cd}
\usepackage[mathscr]{euscript}
\newcommand{\abs}[1]{\left\vert#1\right\vert}
\usepackage[normalem]{ulem}
\usepackage{manfnt}
\usepackage{pgfplots}
\usepackage{pifont}
\usepackage[safe]{tipa}
\usepackage{marvosym}
\usepackage{scalerel,stackengine}
\usepackage{titlesec}
\usepackage{makeidx}
\usepackage{centernot}
\usepackage{xspace}
\usepackage[noauto]{chappg}
\usepackage[OT2,T1]{fontenc}
\usepackage{hyperref}
\usepackage{amsmath}
\usepackage[amsmath,thmmarks]{ntheorem}
\usepackage{tabularx}
\usepackage{yfonts}
\usepackage{epstopdf}
\usepackage{pdfpages}
\newcolumntype{Y}{>{\raggedright\arraybackslash}X}

\newcommand\extgeocze{%
\mathfrak{U}
}

\newcommand\R{%
\mathbb{R}
}

\newcommand\X{%
\mathbb{X}
}
\newcommand\XX{%
\mathbb{X}
}

\newcommand\Y{%
\mathbb{Y}
}

\newcommand\sB{%
\mathcal{B}
}

\newcommand\sH{%
\mathcal{H}
}

\newcommand\sR{%
\mathcal{R}
}

\newcommand\ap{%
\text{ap}
}

\newcommand\Gr{%
\text{Gr}
}

\newcommand\id{%
\text{id}
}

\newcommand\td{%
\text{d}
}
\newcommand\tdx{%
\text{d$x$}
}
\newcommand\tdy{%
\text{d$y$}
}

\newcommand\tA{%
\text{A}
}
\newcommand\tB{%
\text{B}
}

\newcommand\tC{%
\text{C}
}
\newcommand\tE{%
\text{E}
}

\newcommand\Lp{
\text{L}
}

\newcommand\Lm{%
\text{L}
}
\newcommand\tD{%
\text{D}
}

\newcommand\tT{%
\text{T}
}

\newcommand\tV{%
\text{V}
}

\def\thereforex{\boldsymbol{\text{ }
\leavevmode
\lower0.4ex\hbox{\textbullet}
\kern-.9em\raise1.1ex\hbox{\textbullet}
\kern-0.9em\lower0.4ex\hbox{\textbullet}
\hspace{0.1cm}\thinspace\text{ }}}
\def\thereforez{\boldsymbol{\text{ }
\leavevmode
\lower0.4ex\hbox{$\circ$}
\kern-.9em\raise1.1ex\hbox{$\circ$}
\kern-0.9em\lower0.4ex\hbox{$\circ$}
\hspace{0.1cm}\thinspace\text{ }}}
\newcommand\ra{%
\rightarrow
}

\newcommand\ds{%
\displaystyle
}

\newcommand\un[1]{%
\underline{#1}\xspace
}

\newcommand\restr[2]{%
{#1}|{#2}
}

\newcommand\hsx{%
\hspace{0.05cm}
}
\newcommand\hsy{%
\hspace{0.03cm}
}

\newcommand{\norm}[1]{\left\lVert #1 \right\rVert}

\newcommand\BV{%
 \text{BV}
}
\newcommand\BVT{%
 \text{BVT}
}
\newcommand\gBV{%
 \text{gBV}
}
\newcommand\gBVT{%
 \text{gBVT}
}

\newcommand\ov[1]{%
\overline{#1}
}

\newcommand\ovs[1]{
\mkern 1.5mu
\overline{\mkern-2.75mu\mbox{$#1$}\raisebox{3.1mm}{}\mkern-1.5mu}
\mkern 1.5mu
}

\definecolor{ultramarine}{RGB}{0, 32, 96}

\definecolor{darkcerulean}{rgb}{0.3, 0.27, 0.49}

\definecolor{forestgreen}{rgb}{0.0, 0.27, 0.13}
\definecolor{forestgreenweb}{rgb}{0.13, 0.55, 0.13}
\definecolor{deepjunglegreen}{rgb}{0.0, 0.29, 0.29}

\definecolor{midnightblue}{rgb}{0.1, 0.1, 0.44}
\definecolor{midnightgreen}{rgb}{0.0, 0.29, 0.33}
\definecolor{myrtle}{rgb}{0.13, 0.26, 0.12}
\definecolor{darkviolet}{rgb}{0.58, 0.0, 0.83}
\definecolor{darkgreen}{rgb}{0.0, 0.2, 0.13}
\definecolor{officegreen}{rgb}{0.0, 0.5, 0.0}

\definecolor{harvardcrimson}{rgb}{0.79, 0.0, 0.09}
\definecolor{hollywoodcerise}{rgb}{0.96, 0.0, 0.63}
\definecolor{debianred}{rgb}{0.84, 0.04, 0.33}
\definecolor{darkturquoise}{rgb}{0.0, 0.81, 0.82}

\definecolor{darktangernine}{rgb}{1.0, 0.66, 0.07}
\definecolor{aureolin}{rgb}{0.99, 0.93, 0.0}
\definecolor{canaryyellow}{rgb}{1.0, 0.94, 0.0}
\definecolor{amber}{rgb}{1.0, 0.75, 0.0}

\definecolor{urobilin}{rgb}{0.88, 0.68, 0.13}
\definecolor{uscgold}{rgb}{1.0, 0.8, 0.0}

\makeatletter

\def\env@sqcases{%
  \let\@ifnextchar\new@ifnextchar
  \left\lbrack
  \def\arraystretch{1.2}%
  \array{@{}l@{\quad}l@{}}%
}
\makeatother

\usepackage{scalerel,stackengine}
\stackMath
\newcommand\reallywidehat[1]{%
\savestack{\tmpbox}{\stretchto{%
  \scaleto{%
    \scalerel*[\widthof{\ensuremath{#1}}]{\kern-.6pt\bigwedge\kern-.6pt}%
    {\rule[-\textheight/2]{1ex}{\textheight}}
  }{\textheight}%
}{0.5ex}}%
\stackon[1pt]{#1}{\tmpbox}%
}
\makeatletter
\DeclareRobustCommand\widecheck[1]{{\mathpalette\@widecheck{#1}}}
\def\@widecheck#1#2{%
    \setbox\z@\hbox{\m@th$#1#2$}%
    \setbox\tw@\hbox{\m@th$#1%
       \widehat{%
          \vrule\@width\z@\@height\ht\z@
          \vrule\@height\z@\@width\wd\z@}$}%
    \dp\tw@-\ht\z@
    \@tempdima\ht\z@ \advance\@tempdima2\ht\tw@ \divide\@tempdima\thr@@
    \setbox\tw@\hbox{%
       \raise\@tempdima\hbox{\scalebox{1}[-1]{\lower\@tempdima\box
\tw@}}}%
    {\ooalign{\box\tw@ \cr \box\z@}}}
\makeatother

\newtheoremstyle{xx}
  {4pt}
  {0pt}
  {\upshape}
  {\bfseries}
  {}
  { }
  {}
  
\makeatletter 
 \newtheoremstyle{myu}%
  {\upshape\item[ \indent\indent\bf\underline{\theorem@headerfont ##2:}]}%
\makeatother
\makeatletter 
 \newtheoremstyle{myn}%
  {\item[\hskip\labelsep \ \bf ##1 \theorem@headerfont ##2.]}%
\makeatother
\theoremstyle{myn}
\newtheorem{theoremn}{Theorem} 
\theoremstyle{myu}
{\upshape}

\newtheorem{x}[theoremn]{}
\newtheorem{rf}[theoremn]{}

\title{\textbf{Analysis 101:\\
Surfaces and Area}}
\author{Garth Warner\\
Department of Mathematics\\
University of Washington}
\date{}	

\titleformat{\chapter}[display]
{\normalfont\filcenter\huge\bfseries}{}{0pt}{\large}

\titleformat{\chapter}[display]
{\normalfont\filcenter\huge\bfseries}{}{0pt}{\large}

\usepackage[OT2,OT1]{fontenc} 
\newcommand\cyr
{
\renewcommand\rmdefault{wncyr} 
\renewcommand\sfdefault{wncyss} 
\renewcommand\encodingdefault{OT2} 
\normalfont
\selectfont
}

\DeclareTextFontCommand{\textcyr}{\cyr}

\makeindex 
\begin{document}

\maketitle                              

\titlespacing*{\chapter}{0pt}{-50pt}{40pt}
\setlength{\parskip}{0.1em}
\pagenumbering{bychapter}
\setcounter{chapter}{0}
\pagenumbering{bychapter}
\include{_surface_and_area}

\printindex
\end{document}

%% file: _surface_and_area.tex
\begingroup
\fontsize{11pt}{11pt}\selectfont


\[
\textbf{ABSTRACT}
\]
\\

Here one will find a rigorous treatment of the simplest situation in Surface Area Theory, 
viz. 
the nonparametric case with domain the unit square in the plane.  
\\[2cm]

\[
\textbf{ACKNOWLEDGEMENT}
\]

Many thanks to David Clark for his rendering the original transcript into AMS-LaTeX.  
Both of us also thank Judith Clare for her meticulous proofreading.
\newpage

\[
\textbf{SURFACES AND AREA}
\]
\\


\hspace{2.1cm} \ \S X. \quad THE FR\'ECHET PROCESS%
\\[-.1cm]

\hspace{2.1cm} \ \S XX. \hspace{-.0cm} BEFORE THE BEGINNING%
\\[-.1cm]

\hspace{2.1cm} \ \S0. \quad THE BEGINNING%
\\[-.1cm]

\hspace{2.1cm} \ \S1. \quad QUASI LINEAR FUNCTIONS%
\\[-.1cm]

\hspace{2.1cm} \ \S2. \quad LEBESGUE AREA%
\\[-.1cm]

\hspace{2.1cm} \ \S3. \quad GE\"OCZE AREA %
\\[-.1cm]

\hspace{2.1cm} \ \S4. \quad APPROXIMATION THEORY %
\\[-.1cm]

\hspace{2.1cm} \ \S5. \quad TONELLI'S CHARACTERIZATION %
\\[-.1cm]

\hspace{2.1cm} \ \S6. \quad TONELLI'S ESTIMATE%
\\[-.1cm]

\hspace{2.1cm} \ \S7. \quad THE ROLE OF ABSOLUTE CONTINUITY  %
\\[-.1cm]

\hspace{2.1cm} \ \S8. \quad STEINER'S INEQUALITY %
\\[-.1cm]

\hspace{2.1cm} \ \S9. \quad EXTENSION PRINCIPLES %
\\[-.1cm]

\hspace{2.1cm} \S10. \quad ONE VARIABLE REVIEW %
\\[-.1cm]

\hspace{2.1cm} \S11. \quad EXTENDED LEBESGUE AREA   %
\\[-.1cm]

\hspace{2.1cm} \S12. \quad THEORETICAL SUMMARY   %
\\[-.1cm]

\hspace{2.1cm} \S13. \quad VARIANTS   %
\\[-.1cm]

\hspace{3.35cm}   REFERENCES
\\

\[
\]



\endgroup 

\chapter{
$\boldsymbol{\S}$\textbf{X}.\quad THE FR\'ECHET PROCESS}
\setlength\parindent{2em}
\setcounter{theoremn}{0}
\renewcommand{\thepage}{\S X-\arabic{page}}

\qquad
Let $(\X,d)$ be a metric space and let $F:\X \ra [0,+\infty]$ be a lower semicontinuous function.  
Assume:
\\[-.5cm]

\qquad
(A) \ 
For each $x \in \X$, there is a sequence $x_n$ $(n = 1, 2, \ldots)$ in $\X - \{x\}$ converging to x such that

\[
\lim\limits_{n \ra \infty} F(x_n) \ = \ F(x).
\]

Let $(\ov{\X},\bar{d})$ be the completion of $(\X,d)$, the elements $\bar{x}$ 
of which being equivalence classes of Cauchy sequences in $\X$.  
Extend $F$ to a function $\ovs{F}:\ov{\X} \ra [0,+\infty]$ by defining

\[
\ovs{F}(\bar{x}) \ = \ \inf\limits_{\{x_n\} \in \bar{x}} \liminf\limits_{n \ra \infty} F(x_n),
\]
where the infimum is taken over all Cauchy sequences in $\bar{x}$.
\\[-.25cm]

\begin{x}{\small\bf THEOREM} \ 
$\ovs{F}$ is an extension of $F$, i.e., 
\[
\restr{\ovs{F}}{\X} \ = \ F.
\]
Moreover $\ovs{F}$ is lower semicontinuous and in addtion is unique.
\\[-.25cm]
\end{x}

\begin{x}{\small\bf \un{N.B.}} \ 
$\ovs{F}$ has the following property: 
\\[-.25cm]

\qquad 
(B) \quad For each $\bar{x} \in \ov{\X}$, there is a Cauchy sequence $\{x_n\} \in \bar{x}$ such that

\[
\qquad\qquad \lim\limits_{n \ra \infty} F(x_n) \ = \ \ovs{F}(\bar{x}).
\]
\\[-1.5cm]
\end{x}

To recapitulate:
\\[-.25cm]

\begin{x}{\small\bf SCHOLIUM} \ 
Every nonnegative, extended real valued, lower semicontinuous function on a metric space $\X$ with property 
(A) can be extended to a unique 
lower semicontinuous function on the completion $\ov{\X}$ of $\X$ with property (B).
\\[-.25cm]
\end{x}

\begin{x}{\small\bf EXAMPLE} \ 
Consider

\[
\begin{cases}
\ \X \ = \ \ ]0,1[ \qquad (d(x,y) = \abs{x - y})
\\[4pt]
\ \ov{\X} \ = \  [0,1] \qquad (\bar{d}(\bar{x},\bar{y}) = \abs{\bar{x} - \bar{y}})
\end{cases}
\]
and

\[
\begin{cases}
\ F \ = \ \id_{\X}
\\[4pt]
\ \ovs{F} \ = \ \id_{\ov{\X}}
\end{cases}
.
\]

\end{x}

\chapter{
$\boldsymbol{\S}$\textbf{XX}.\quad BEFORE THE BEGINNING}
\setlength\parindent{2em}
\setcounter{theoremn}{0}
\renewcommand{\thepage}{\S XX-\arabic{page}}

\ 
\\[.25cm]
\qquad 
\textbf{Definition} \ (J. A. Serret, \textit{Cours de Calcul Differential et Integral}  Vol II, (1868) p. 296.) 
Soit une portion de surface courbe termin\'ee par un contour $C$; \ 
nous nommerons aire de cette surface la limite $S$ vers laquelle tend l'aire d'une surface polydrale inscrite form\'ee 
de faces triangulaires et termin\'ee par un contour polygonal $\Gamma$ ayant pour limite le contour $C$.
\\[-.25cm]

\textit{
Area of a surface $S$ bounded by a curve $C$; \ 
The area is the limit of the elementary areas of the inscribed polyhedral surfaces $P$ bounded by a curve $\Gamma$ as 
$P \ra S$ and $\Gamma \ra C$, where this limit exists and is independent of the particular sequence of inscribed polygons.  
}
\\

The above `definition' is the obvious generalization of the standard definition of arc length of a curve to the area of a surface.  
However, the nature of the area of a surface is more subtle (and thus more interesting) than this simple definition allows.  
The following example, due to H. A. Schwarz (and independently by G. Peano), shows Serret's definition of the area of a surface 
to be untenable. 
Lore (apochryphal or not) is that Schwarz realized this counterexample while imagining polyhedral approximations to a Christmas tree ornament on Chrismas eve.  
The result was described by Schwarz in a letter dated 25-26 December, 1880,  to A. Genoccchi.  
\\

\begin{spacing}{1.75}
\textbf{Example} \ (H. A. Schwarz)

Partition the cylinder of height 1 and radius one, 
vertically, into $m$ slices each of height $\ds\frac{1}{m}$, and   
horizontally, into n sectors each of angle $\ds\frac{2 \pi}{n}$.

(There is no generality to be gained in assuming the cylinder has height $h$ and radius $r$.) 
\end{spacing} 
\newpage

\vspace{-1.5cm}
\begin{figure}[h]
\centerline{\includegraphics[width=10cm, height=7cm]{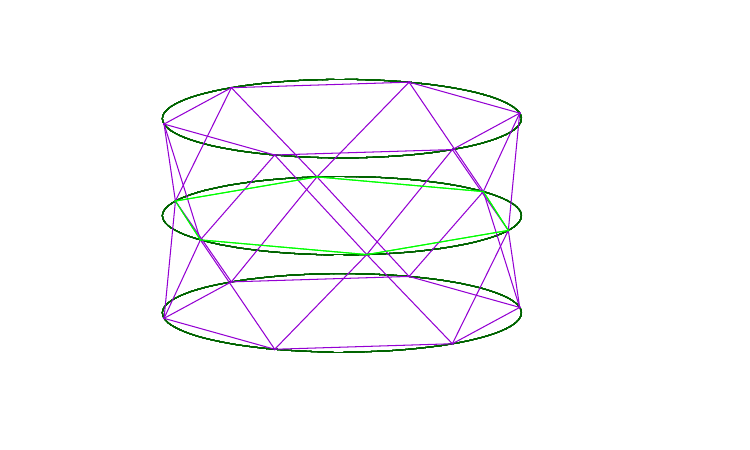}}
\end{figure}

\vspace{-.5cm}
Let $P(m,n)$ be the inscribed polyhedral surface consisting of $2 \hsy m \hsy n$ inscribed triangles per the partition described above.  
Let $\{T(m, n)_j : j = 1, \ldots, 2 \hsy m \hsy n\}$ be the set of triangles comprising the polyhedal surface $P(m,n)$.  
Each of the  triangles $T(m, n)_j $ has the same area.  
Denote this common area by 
\text{Area$(T(m, n))$}.
\\[-.25cm]

\qquad $\implies$
\\[-1.25cm]
\begin{align*}
\text{Area $(P(m,n))$} \ 
&=\ 
2 \hsy m \hsy n \cdot \text{Area$(T(m, n))$}
\end{align*}
\\[-2.5cm]

{
\begin{center}
\includegraphics[width=14cm, height=8cm]{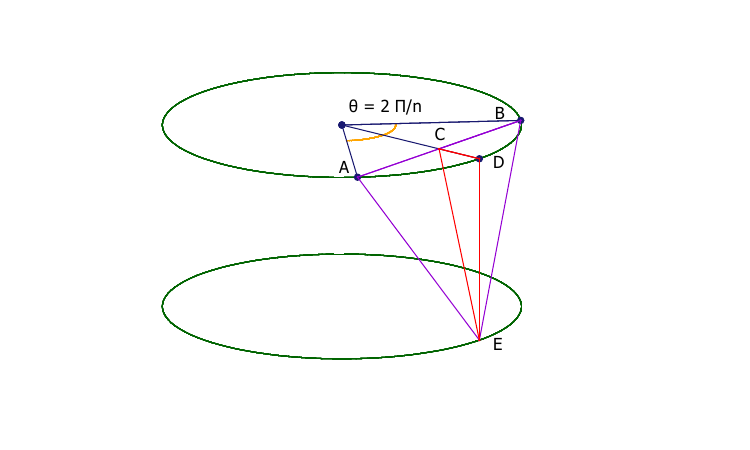}
\end{center}
}

\vspace{-.75cm}
\un{Compute}: \ 
\allowdisplaybreaks
\begin{align*}
\text{Area$(T(m, n))$} \ 
&= \  
\text{Area ($\triangle$ (A,B,E))} \ 
\\[18pt]
&= \ 
\frac{1}{2} \ \norm{\tA - \tB} \cdot \norm{\tC - \tE}
\\[18pt]
&= \
\frac{1}{2} \ 2 \sin\bigg(\frac{\pi}{n}\bigg) \ 
\hsx \cdot  \hsx 
\sqrt{\norm{\tC - \tD}^2 + \norm{\tD - \tE}^2)}
\\[18pt]
&= \
\sin\bigg(\frac{\pi}{n}\bigg) \ 
\hsx \cdot  \hsx 
\sqrt{\bigg(1 - \cos\bigg(\frac{\pi}{n}\bigg)\bigg)^ 2 + \bigg(\frac{1}{m}\bigg)^2}
\\[18pt]
&\approx \
\bigg(\frac{\pi}{n}\bigg) 
\hsx \cdot  \hsx 
\sqrt{\bigg(1 - \bigg(1 - \frac{1}{2} \hsy \bigg(\frac{\pi}{n}\bigg)^2\bigg)^2  \hsx + \hsx \bigg(\frac{1}{m}\bigg)^2}
\\[18pt]
&= \
\bigg(\frac{\pi}{n}\bigg) 
\hsx \cdot  \hsx 
\sqrt{\frac{1}{2} \hsx \frac{\pi^4}{n^4}  \hsx + \hsx \frac{1}{m^2}}.
\end{align*}

Therefore
\begin{align*}
\text{Area $(P(m,n))$} \ 
&=\ 
2 \hsy m \hsy n 
\cdot 
\text{Area$(T(m, n))$}
\\[18pt]
&=\ 
2 \hsy m \hsy n \cdot 
\text{Area $\triangle$ (A,B,E)} \ 
\\[18pt]
&\approx \
2 \hsy m \hsy n \cdot 
\bigg(\frac{\pi}{n}\bigg) 
\hsx \cdot  \hsx 
\sqrt{\frac{1}{2} \hsx \frac{\pi^4}{n^4}  \hsx + \hsx \frac{1}{m^2}}
\\[18pt]
&\approx \
2 \hsy \pi 
\hsx \cdot  \hsx 
\sqrt{m^2 \hsy \frac{1}{2} \hsx \frac{\pi^4}{n^4}  \hsx + \hsx 1}.
\end{align*}

\un{Case}: \ Let $n \ra \infty$ first and then let $m \ra \infty$:
\begin{align*}
\lim\limits_{m \ra \infty} \ \lim\limits_{n \ra \infty} 
\text{Area $(P(m,n))$} \ 
&=\
\lim\limits_{m \ra \infty} \ \lim\limits_{n \ra \infty} 
2 \hsy \pi 
\hsx \cdot  \hsx 
\sqrt{m^2 \hsy \frac{1}{2} \hsx \frac{\pi^4}{n^4}  \hsx + \hsx 1}
\\[15pt]
&= \ 
2 \hsy \pi.
\end{align*}

\un{Case}: \ Let $m \ra \infty$ first and then let $n \ra \infty$:
\begin{align*}
\lim\limits_{n \ra \infty} \ \lim\limits_{n \ra \infty} 
\text{Area $(P(m,n))$} \ 
&=\
\lim\limits_{n \ra \infty} \ \lim\limits_{n \ra \infty} 
2 \hsy \pi 
\hsx \cdot  \hsx 
\sqrt{m^2 \hsy \frac{1}{2} \hsx \frac{\pi^4}{n^4}  \hsx + \hsx 1}
\\[15pt]
&= \ 
\infty.
\end{align*}

\un{Case}: \ Let $n = m \ra \infty$:
\begin{align*}
\lim\limits_{n = m \ra \infty} \ 
\text{Area $(P(n,n))$} \ 
&=\
\lim\limits_{n = m \ra \infty} \ 
2 \hsy \pi 
\hsx \cdot  \hsx 
\sqrt{\frac{1}{2} \hsx \frac{\pi^4}{n^2}  \hsx + \hsx 1}
\\[15pt]
&= \ 
2 \hsy \pi.
\end{align*}

\un{Case}: \ Let $m = c \hsx \ds\frac{\sqrt{2}}{\pi^2} \hsx n^2 \ra \infty$ first and then let $m \ra \infty$ \ $(c \in \R_{\geq 0})$.
\begin{align*}
\lim\limits_{m = c n^2 \ra \infty} \ 
\text{Area $(P(n,n))$} \ 
&=\
\lim\limits_{n = m \ra \infty} \ 
2 \hsy \pi 
\hsx \cdot  \hsx 
\sqrt{c^2  \hsx + \hsx 1}
\\[15pt]
&= \ 
2 \hsy \pi.
\end{align*}

\un{Moral}: \ 
Depending on how $m$, $n \ra \infty$;
\[
\lim\limits_{m, \hsy n \ra \infty} \ \text{Area $(P(m,n))$}
\]
can assume any value $c \in [2 \hsy \pi, +\infty]$.  
\\[-.5cm]

In particular, the limit of the areas of inscribed polygons is not unique, thus rendering Serret's definition of surface area untenable.  
However, contained in this example is a hint towards the more satisfactory definition later set forth by Lebesgue.




\chapter{
$\boldsymbol{\S}$\textbf{0}.\quad THE BEGINNING}
\setlength\parindent{2em}
\setcounter{theoremn}{0}
\renewcommand{\thepage}{\S0-\arabic{page}}

\qquad 
Traditionally, a \un{$k$-surface in $n$-space} $(k \leq n)$ is an ordered pair $S = (A, \un{f})$, 
where $A$ is a subset of $\R^k$ with a nonempty interior (subject to certain restrictions) and 
$\un{f}$ is a function from $A$ to $\R^n$, 
i.e., 
$\un{f} : A \ra \R^n$, thus

\[
\un{f}
\ = \ 
(f_1, \ldots, f_n).
\]
\\[-1.25cm]

\begin{x}{\small\bf \un{N.B.}} \ 
If $k = n$, then $\un{f}$ is said to be \un{flat}.
\\[-.25cm]
\end{x}

\begin{x}{\small\bf REMARK} \ 
If $k = 1$ and $A = [a,b]$, then $\un{f}$ is just a curve.
\\[-.25cm]
\end{x}

In this account, we shall take $k = 2$ and $n = 3$, thus

\[
\un{f} \hsx : \ 
\begin{cases}
\ f_1 : A \ra \R
\\[4pt]
\ f_2 : A \ra \R
\\[4pt]
\ f_3 : A \ra \R
\end{cases}
.
\]
\\[-.75cm]

\begin{x}{\small\bf \un{N.B.}} \ 
There are associated flat maps, viz.

\[
\begin{cases}
\ x = 0, \hspace{1.1cm} \quad y = f_2 (u, v), \quad z = f_3 (u, v)
\\[4pt]
\ x = f_1 (u, v), \quad y = 0, \hspace{1.1cm}  \quad z = f_3 (u, v)
\\[4pt]
\ x = f_1 (u, v), \quad y = f_2 (u, v), \quad z = 0
\end{cases}
,
\]
where $(u, v) \in A$.
\\[-.25cm]
\end{x}

In what follows, we do not intend to opererate ``in general'' but instead will specialize matters to the so-called 
``nonparametric'' situation. 
\\[-.25cm]


Put 

\[
Q 
\ = \ 
[0,1] \hsx \times \hsx [0,1] \subset \R^2 
\qquad 
(0 \leq x \leq 1, 0 \leq y \leq 1).
\]
\\[-1cm]

\begin{x}{\small\bf DEFINITION} \ 
A \un{nonparametric} 2-surface in 3-space is an ordered pair $S_f = (Q, \un{f})$, where

\[
\un{f} (x, y) = (x, y, f(x, y)), \quad f: Q \ra \R
\]
is a function, thus

\[
\begin{cases}
\ 
f_1(x, y) = x
\\[8pt]
\ 
f_2(x, y) = y
\end{cases}
\
f_3 (x, y) = f (x, y).
\]
\\[-1.cm]
\end{x}

\begin{x}{\small\bf REMARK} \ 
Every function $f: Q \ra \R$ determines a nonparametric surface $S_f$.  
Because of this, the focus is on $f$, not $S_f$. 
\\[-.25cm]
\end{x}

Restricting matters to $Q$ more or less eliminates the topological apsects of the theory, 
thus the discussion is ``pure analysis'', 
there being two aspects to the development, viz. 

\[
\begin{cases}
\ 
\text{PART 1: \ The Continuous Case, $f \in C(Q)$.}
\\[4pt]
\ 
\text{PART 2: \ The Integrable Case, $f \in \Lp^1(Q)$.}
\end{cases}
\]
\\[-.75cm]

\begin{x}{\small\bf EXAMPLE} \ 
Define $f: Q \ra \R$ by the prescription

\[
\begin{cases}
\ 
0 \hspace{0.75cm} \Big(0 \leq x \leq \ds \frac{1}{2}\hsx \Big)
\\[18pt]
\ 
1 \hspace{0.75cm} \Big(\ds \frac{1}{2} < x \leq 1 \hsx\Big)
\end{cases}
.
\]
Then $f$ is not continuous but it is integrable.
\end{x}

\chapter{
$\boldsymbol{\S}$\textbf{1}.\quad QUASI LINEAR FUNCTIONS}
\setlength\parindent{2em}
\setcounter{theoremn}{0}
\renewcommand{\thepage}{\S1-\arabic{page}}

\begin{x}{\small\bf DEFINITION} \ 
A \underline{quasi linear function} is a continuous function 
$\Pi:Q \rightarrow \R$ 
for which there exists a decomposition $D$ of $Q$ into a finite number of nonoverlapping triangles \hsx 
$T_1, T_2, \ldots, T_n$ such that $\Pi$ is linear in each of these triangles, thus
\[
\Pi(x,y) = a_i x + b_i y + c_i \qquad ((x,y) \in T_i),
\]
the $a_i, b_i, c_i$ being real numbers.
\\[-.5cm]
\end{x}

\begin{x}{\small\bf EXAMPLE} \ 
A constant function

\[
f(x,y) 
\ = \  
C  \qquad ((x,y) \in Q)
\]
is quasi linear.
\\[-.5cm]
\end{x}

Suppose that $\Pi:Q \rightarrow \R$ is quasi linear $-$then $\Pi$ maps each $T_i$ 
into a triangle $\Delta_i \subset \R^3$ (possibly a segment or a point).
\\[-.25cm]

\begin{x}{\small\bf NOTATION} \ 
Let $\abs{\Delta_i}$ stand for the area of $\Delta_i$.
\\[-.25cm]
\end{x}

\begin{x}{\small\bf DEFINITION} \ 
The \underline{elementary area} of a quasi linear function $\Pi:Q \rightarrow \R$ is the sum

\[
a(\Pi) 
\ \equiv \ 
\sum \ \abs{\Delta_i},
\]
where $\ds \sum$ is taken over the $T_i \in D$.
\\[-.5cm]
\end{x}

\begin{x}{\small\bf NOTATION} \ 
Let $\abs{T_i}$ stand for the area of $T_i$.
\\[-.5cm]
\end{x}

\begin{x}{\small\bf \un{N.B.}} \ 
Let
\[
(u_1, v_1), (u_2, v_2), (u_3, v_3)
\]
be the vertices of $T_i$ in Q $-$then


\[\abs{T_i} \ = \  \frac{1}{2} \hsx  \left| \hsx \det \ 
\begin{bmatrix}
\hsx u_1 &v_1 &1 \hsx\\
\hsx u_2 &v_2 &1 \hsx\\
\hsx u_3 &v_3 &1 \hsx\\
\end{bmatrix}
\right|.
\]
\\[-1cm]
\end{x}

\begin{x}{\small\bf LEMMA} \ 

\[
\abs{\Delta_i} = \abs{T_i} (1 + a_i^2 + b_i^2)^{1/2}.
\]
Therefore

\[
a(\Pi) = \sum\limits_i \abs{T_i} (1 + a_i^2 + b_i^2)^{1/2}.
\]
\\[-1cm]
\end{x}

\begin{x}{\small\bf SCHOLIUM} \ 

\[
a(\Pi) 
\ = \ 
\iint\limits_Q \ [1 + (\partial \Pi/\partial x)^2 + (\partial \Pi /  \partial y)^2]^{1/2} \ \td x \hsy \td y.
\]
\\[-1cm]
\end{x}

It follows from this that $a(\Pi)$ is independent of the subdivision $D$ of $Q$ into triangles of linearity for $\Pi$.
\\[-.25cm]

\begin{x}{\small\bf REMARK} \ 
A quasi linear function $\Pi:Q \rightarrow \R$ is Lipschitz continuous and 

\[
\sH^2 (\Gr_\Pi (Q) )
\ = \ 
\iint\limits_Q \ [1 + (\partial \Pi/\partial x)^2 + (\partial \Pi /  \partial y)^2]^{1/2} \ \td x \hsy \td y.
\]
\\[-1cm]
\end{x}

\begin{x}{\small\bf LEMMA} \ 
Per uniform convergence, the elementary area is lower semicontinuous on the set of quasi linear functions. 
\end{x}

\chapter{
$\boldsymbol{\S}$\textbf{2}.\quad LEBESGUE AREA}
\setlength\parindent{2em}
\setcounter{theoremn}{0}
\renewcommand{\thepage}{\S2-\arabic{page}}

\qquad
Recall that
\[
Q \ = \ [0,1] \times [0,1] \ \subset \ \R^2 \qquad (0 \leq x \leq 1, 0 \leq y \leq 1).
\]
\\[-1.25cm]

\begin{x}{\small\bf LEMMA} \ 
Let $F:Q \ra \R$ be a continuous function $-$then there exists a sequence
\[
\xi \ = \ \{\Pi_n: n = 1, 2, \ldots\}
\]
of quasi linear functions $\Pi_n:Q \ra \R$ such that $\Pi_n \ra f$ uniformly $(n \ra \infty)$.
\\[-.25cm]
\end{x}

\begin{x}{\small\bf NOTATION} \ 
Given a continuous function $f:Q \ra \R$, denote by $\Xi$ the collection of all sequences
\[
\xi \ = \ \{\Pi_n: n = 1, 2, \ldots\}
\]
of quasi linear functions $\Pi_n:Q \ra \R$ such that $\Pi_n \ra f$ uniformly $(n \ra \infty)$.
\\[-.25cm]
\end{x}

\begin{x}{\small\bf \un{N.B.}} \ 
The preceeding lemma ensures that $\Xi$ is nonempty.
\\[-.25cm]
\end{x}

\begin{x}{\small\bf DEFINITION} \ 
The \underline{Lebesgue area} $L_Q[f]$ of a continuous function $f:Q \ra \R$ is the entity
\[
\inf\limits_{\xi  \in \hsy \Xi} \ \liminf\limits_{n \ra \infty} \ a(\Pi_n).
\]
\\[-1.cm]
\end{x}

\begin{x}{\small\bf REMARK} \ 
This definition and the considerations that follow are an instance of the Fr\'echet process:  
Take for $\X$ the quasi linear functions on $Q$, take for $d$ the metric defined by the prescription

\[
d(\Pi_1, \Pi_2) \ = \ \sup \abs{\Pi_1(x,y) - \Pi_2(x,y)},
\]
and take for $F$ the elementary area $-$then the completion $\ov{\X}$ of $\X$ is $C(Q)$, the 
set of continuous functions on $Q$, and the extension $\ovs{F}$ of $F$ assigns to each $f \in C(Q)$ its Lebesgue area:
\[
\qquad \bar{F}(f) \ = \ L_Q[f].
\]
\\[-1cm]
\end{x}

\begin{x}{\small\bf CONSISTENCY PRINCIPLE}\ 
The elementary area of a quasi linear function $\Pi:Q \ra \R$ equals its Lebesgue area.
\\[-.25cm]
\end{x}

\begin{x}{\small\bf LEMMA} \ 
There is at least one $\xi \in \Xi$ such that
\[
a(\Pi_n) \ra L_Q[f] \qquad (n \ra \infty).
\]

PROOF \ 
There are two possibilities:
\[
\begin{cases}
\ L_Q[f] < +\infty
\\
\quad \text{or} 
\\
\ L_Q[f] = +\infty
\end{cases}
.
\]
Matters are manifest if $L_Q[f] = +\infty$, so assume that $L_Q[f] < +\infty$.  
Given any positive integer $n$, there exists a sequence $\{\Pi_m: m = 1, 2, \ldots\}$ 
such that for $m \ra \infty$, $\Pi_m \ra f$ uniformly and

\[
\liminf\limits_{m \ra \infty} a(\Pi_m) \ < \ L_Q[f] + \frac{1}{n},
\]
thus there is an $m$ such that

\[
\norm{\Pi_m - f}_\infty \ < \ \frac{1}{n}
\]
and
\[
a(\Pi_m) \ < \ \ L_Q[f] + \frac{1}{n}.
\]
This $m$ depends on $n$.  
Write $\Pi(n)$ in place of $\Pi_m$ $-$then

\[
\norm{\Pi(n) - f}_\infty \ < \ \frac{1}{n}
\]
and

\[
a(\Pi(n)) \ < \ \ L_Q[f] + \frac{1}{n}.
\]
Let now $n \ra \infty$ to conclude that

\[
\Pi(n) \ra f
\]
uniformly and

\[
\limsup\limits_{n \ra \infty} a(\Pi(n)) \ \leq \ L_Q[f].
\]
On the other hand, 

\[
\ L_Q[f] \ \leq \ \liminf\limits_{m \ra \infty} a(\Pi(n)).
\]
Hence the lemma.
\\[-.25cm]
\end{x}

\begin{x}{\small\bf \un{N.B.}} \ 
This result is known as the proper sequential limit principle.
\\[-.25cm]
\end{x}

\begin{x}{\small\bf THEOREM} \ 
Let $f:Q \ra \R$ be a continuous function.  
Suppose that 
$f_n:Q \ra \R$ $(n = 1, 2, \ldots)$ is a sequence of continuous functions such that $f_n \ra f$ uniformly $-$then

\[
\ L_Q[f] \ \leq \ \liminf\limits_{n \ra \infty} L[f_n].
\]

PROOF \ 
Assume without loss of generality that

\[
\begin{cases}
\ 
\liminf\limits_{n \ra \infty} L_Q[f_n] \ < \ +\infty
\\[4pt]
\quad \text{and} 
\\[4pt]
\ 
 L_Q[f_n] \hspace{1.35cm}< +\infty
\end{cases}
\qquad  (\forall \ n).
\]
Given $n$, choose per supra a sequence 
$\{\Pi_{n \hsy m}:m = 1, 2, \ldots\}$ 
of quasi linear functions uniformly convergent to 
$f_n$ $(m \ra \infty)$ 
with

\[
a(\Pi_{n \hsy m}) \ra L_Q[f_n] \qquad (m \ra \infty).
\]
Accordingly

\[
\delta_{n \hsy m} \ \equiv \ \norm{\Pi_{n \hsy m} - f_n}_\infty \ra 0 \qquad (m \ra \infty)
\]
and for each $n$ there exists an integer $m = m(n)$ such that

\[
\delta_{n \hsy m} \ < \ \frac{1}{n} \qquad \text{and} \qquad \abs{a(\Pi_{n \hsy m}) - L_Q[f_n]} \ < \ \frac{1}{n}.
\]
Next, $\forall \ w \in Q$,
\allowdisplaybreaks
\begin{align*}
\abs{\Pi_{n \hsy m}(w) - f(w)} \ 
&\leq \ 
\norm{\Pi_{n \hsy m} - f_n}_\infty + \norm{f_n - f}_\infty
\\[11pt]
&\leq \ 
\delta_{n \hsy m} + \norm{f_n - f}_\infty
\\[11pt]
&< \ 
\frac{1}{n} + \norm{f_n - f}_\infty
\\[11pt]
&
\ra 0 \qquad (n \ra \infty).
\end{align*}
Put

\[
\Pi_n^\prime \ = \ \Pi_{n \hsy m}
\]
and let

\[
\xi^\prime \ = \ \{\Pi_n^\prime: n = 1, 2, \ldots\},
\]
so $\xi^\prime \in \Xi$.  And
\allowdisplaybreaks
\begin{align*}
L_Q[f] \ 
&\leq \ 
\liminf\limits_{n \ra \infty} \ a(\Pi_n^\prime) 
\\[11pt]
&= \ 
\liminf\limits_{n \ra \infty} \ (a(\Pi_n^\prime) - L_Q[f_n] + L_Q[f_n])
\\[11pt]
&= \ 
\lim\limits_{n \ra \infty}  \ (a(\Pi_n^\prime) - L_Q[f_n]) + \liminf\limits_{n \ra \infty} \ L_Q[f_n]
\\[11pt]
&= \ 
0 + \liminf\limits_{n \ra \infty} \  L_Q[f_n]
\\[11pt]
&= \ 
\liminf\limits_{n \ra \infty} \ L_Q[f_n].
\end{align*}

Therefore Lebesgue area is a lower semicontinuous functional in the class of continuous functions 
(the underlying convergence being uniform).
\\[-.5cm]

[Note: \  
It can be shown that Lebesgue area is a lower semicontinuous functional in the class of continuous functions relative to pointwise convergence.]
\\[-.25cm]
\end{x}

Here is a simple application: \  If $\forall \ n$, $L_Q[f_n] \ \leq \ L_Q[f]$, then \  $L_Q[f_n] \ra L_Q[f]$.  
In fact, 
\[
\limsup\limits_{n \ra \infty} L_Q[f_n] \ \leq \ L_Q[f]
\]
while on the other hand,
\[
\liminf\limits_{n \ra \infty} L_Q[f_n] \ \geq \ L_Q[f].
\]
\\[-1cm]

\begin{x}{\small\bf LEMMA} \ 
Let $L^*$ be a functional in the class of continuous functions which is lower semicontinuous 
per uniform convergence and has the property that for every quasi linear $\Pi$,
\[
L^*[\Pi] \ = \ a(\Pi).
\]
Then for every $f$,

\[
L^*[f] \ \leq \ L_Q[f].
\]

PROOF \ 
Choose $\xi \in \Xi$ such that

\[
a(\Pi_n) \ra  L_Q[f] \qquad (n \ra \infty)
\]
and note that
\begin{align*}
L^*[f] \ 
&\leq \ 
\liminf\limits_{n \ra \infty} \ L^*[\Pi_n]
\\[11pt]
&= \ 
\liminf\limits_{n \ra \infty} \ a(\Pi_n)
\\[11pt]
&\leq \  L_Q[f].
\end{align*}
\end{x}

\chapter{
$\boldsymbol{\S}$\textbf{3}.\quad GE\"OCZE AREA}
\setlength\parindent{2em}
\setcounter{theoremn}{0}
\renewcommand{\thepage}{\S3-\arabic{page}}

\qquad
The setting for the notion of Lebesgue area is the unit square
\[
Q \ = \ [0,1] \times [0,1].
\]
However there is no difficulty in extending matters to oriented rectangles $R \subset Q:$

\[
\begin{cases}
\ a \leq x \leq b \qquad (a < b)
\\[8pt]
\ c \leq y \leq d \qquad (c < d)
\end{cases}
, \ \abs{R} = (b - a) (d - c).
\]
\\[-1cm]

The theory thus formulated applies to any real valued continuous function on $R$.  
In particular:  
Given a continuous function $f:Q \ra \R$ let $f_R$ be its restriction to $R$ and denote its Lebesgue area per $R$ by the symbol $L_Q[f_R]$.
\\[-.25cm]

Introduce

\[
\begin{cases}
\ \ds
G_\X(f;R) 
\ = \ 
\int\limits_a^b \ \abs{f(x,d) - f(x,c)} \ \tdx 
\\[26pt]
\ \ds
G_\Y(f;R) 
\ = \ 
\int\limits_c^d \ \abs{f(b,y) - f(a,y)} \ \tdy 
\end{cases}
\]
and put

\[
\Gamma(f;R) 
\ = \ 
[(G_\X(f;R) )^2 + (G_\Y(f;R) )^2 + \abs{R}^2]^{1/2}.
\]
\\[-1cm]

\begin{x}{\small\bf LEMMA} \ 
\[
\Gamma(f;R) \ \leq \ L_Q[f_R].
\]
\end{x}

Let $D$ be a subdivision of $Q$ into nonoverlapping oriented rectangles $R$ (lines parallel to the coordinate axes).
\\[-.25cm]

\begin{x}{\small\bf DEFINITION} \ 
The \un{sum of Ge\"ocze} is the expression
\[
G(f;D) \ = \ \sum \ \Gamma(f;R),
\]
the summation being taken over the rectangles $R$ in $D$.
\\[-.5cm]
\end{x}

So

\[
G(f;D) \ \leq \ \sum \ L_Q[f_R].
\]
And

\[
\sum\  L_Q[f_R] \ \leq \ L_Q[f].
\]
Therefore

\[
G(f;D) \ \leq \ L_Q[f].
\]
\\[-1cm]

\begin{x}{\small\bf NOTATION} \ 
Put

\[
\Gamma_Q[f] \ = \ \sup\limits_{D} \ G(f;D),
\]
the \un{Ge\"ocze area} of $f$.
\\[-.25cm]

Then $\forall$ D,
\[
G(f;D) \ \leq \ L_Q[f]
\]

\qquad
$\implies$
\[
\Gamma_Q[f] \ \leq \ L_Q[f].
\]
\\[-1.25cm]

[Note: \  
This inequality is trivial if $L_Q[f] = +\infty$, thus there is no loss in generality in assuming that $L_Q[f] < +\infty.]$
\\[-.25cm]
\end{x}

\begin{x}{\small\bf THEOREM} \ 
\[
\Gamma_Q[f] \ = \ L_Q[f].
\]
\end{x}

This assertion is nontrivial, the first step being to establish it when
\[
\frac{\partial f}{\partial x} = p(x,y), \quad \frac{\partial f}{\partial y} = q(x,y)
\]
exist in $Q$ and are continuous there.
\\[-.25cm]

\textbullet \quad
Write
\begin{align*}
G_\X(f;R) \ 
&= \ 
\int\limits_a^b \abs{f(x,d) - f(x,c)} \ \tdx 
\\[15pt]
&= \ 
(b - a) \abs{f(\xi,d) - f(\xi,c)} \qquad (a \leq \xi \leq b)
\\[15pt]
&= \ 
(b - a) (d - c) \abs{q(\xi,\eta)} \qquad\quad (c \leq \eta \leq d)
\\[15pt]
&= \ 
\abs{R} \abs{q(\xi,\eta)}.
\end{align*}
\\[-1.25cm]

\textbullet \quad 
Write
\begin{align*}
G_\Y(f;R) \ 
&= \ 
\int\limits_c^d \abs{f(b,y) - f(a,y)} \ \tdy 
\\[15pt]
&= \ (d - c) \abs{f(b, \mu) - f(a, \mu)} \qquad (c \leq \mu \leq d)
\\[15pt]
&= \ 
(d - c) (b - a) \abs{p(\nu, \mu)} \qquad\quad\  (a \leq \nu \leq b)
\\[15pt]
&= \ 
\abs{R} \abs{p(\nu, \mu)}.
\end{align*}
Consequently
\begin{align*}
\Gamma(f;R) \ 
&= \ 
[1 + p(\nu, \mu)^2 + q(\xi,\eta)^2]^{1/2} \abs{R} 
\\[11pt]
&= \ 
[1 + p(\xi,\eta)^2 + q(\nu, \mu)^2]^{1/2} \abs{R} + \varepsilon_R \abs{R}, 
\end{align*}
where $\varepsilon_R$ tends to zero with the diameter of $R$.
\\[-.5cm]

Let again $D$ be a subdivision of $Q$ into nonoverlapping oriented rectangles $R$ $($lines parallel to the coordinate axes$)$.  
Since $\sum \abs{R} = \abs{Q} = 1$, it follows that
\allowdisplaybreaks
\begin{align*}
G(f;D) \ 
&= \ 
\sum \ \Gamma(f;R)
\\[15pt]
&= \ 
\sum \ [1 + p(\xi,\eta)^2 + q(\nu, \mu)^2]^{1/2} \abs{R} + \varepsilon.
\end{align*}
Here $\varepsilon \ra 0$ when $\delta \ra 0$ $(\delta$ being the maximum diameter of the rectangle $R$ in $D$).
\\[-.5cm]

Replace now $D$ by a sequence $\{D_n\}$ and assume that $\delta_n \ra 0 \ (n \ra \infty)$ $-$then the sum

\[
\sum \ [1 + p(\xi,\eta)^2 + q(\nu, \mu)^2]^{1/2} \abs{R}
\]
tends to the integral

\[
\iint\limits_Q \ (1 + p^2 + q^2)^{1/2} \ \tdx \hsy \tdy ,
\]
hence

\[
\lim\limits_{n \ra \infty} \ G(f;D_n) 
\ = \ 
\iint\limits_Q \ (1 + p^2 + q^2)^{1/2} \ \tdx \hsy \tdy 
\]
or still,
\begin{align*}
\Gamma_Q[f] \ 
&\geq \ 
\iint\limits_Q \ (1 + p^2 + q^2)^{1/2} \ \tdx \hsy \tdy 
\\[15pt]
&\equiv \ 
L_Q[f] \qquad (\text{see below}).
\end{align*}
But, as has been noted above, it is always the case that

\[
\Gamma_Q[f] \ \leq \ L_Q[f].
\]
So in the end,
\[
\Gamma_Q[f] \ = \ L_Q[f].
\]
\\[-1cm]

\begin{x}{\small\bf CONSTRUCTION} \ 
There is a $\xi \in \Xi$ such that
\\[-.5cm]

\[
a(\Pi_n) \quad (n \ra \infty)\  \longrightarrow  \iint\limits_Q \ (1 + p^2 + q^2)^{1/2} \  \tdx \hsy \tdy .
\]
\\[-1cm]
\end{x}

\begin{x}{\small\bf LEMMA} \ 
\[
L_Q[f] \ = \ \iint\limits_Q \ (1 + p^2 + q^2)^{1/2} \ \tdx \hsy \tdy .
\]

PROOF \ 
\begin{align*}
\iint\limits_Q \ (1 + p^2 + q^2)^{1/2} \ \td x \hsy \td y \ 
&\leq \ 
\Gamma_Q[f]
\\
&\leq \ 
L_Q[f]
\\[15pt]
&\leq \ 
\liminf\limits_{n \ra \infty} \ a(\Pi_n)
\\[15pt]
& = \ 
\lim\limits_{n \ra \infty} \ a(\Pi_n)
\\[15pt]
&= \ \iint\limits_Q \ (1 + p^2 + q^2)^{1/2} \  \ \tdx \hsy \tdy .
\end{align*}
\\[-1cm]
\end{x}


\begin{x}{\small\bf EXAMPLE} \ 
Suppose that $f(x,y)$ is independent of y $-$then 
$\ds\frac{\partial f}{\partial y} = 0$  
and 
$\ds\frac{\partial f}{\partial x} = f^\prime(x)$, hence

\[
\iint\limits_Q \ (1 + p^2 + q^2)^{1/2} \ \tdx \hsy \tdy  
\ = \ 
\int\limits_0^1 \ (1 + (f^\prime(x))^2)^{1/2} \ \tdx .
\]
\\[-1cm]
\end{x}

It remains to establish that
\[
\Gamma_Q[f] \ = \ L_Q[f]
\]
in general.  
To this end, denote by \un{$Q$} a concentric square completely contained in the interior of $Q$, let $0 < h < \frac{1}{2}$, put

\[
Q_h: \ 
\begin{cases}
\ h \leq x \leq 1 - h
\\[4pt]
\ h \leq y \leq 1 - h
\end{cases}
,
\]
and assume that for $h$ sufficiently small, $\un{Q} \subset Q_h$ 
$-$then there exists a continuous function $f_h:Q_h \ra \R$ with the following properties.
\begin{align*}
\qquad 
&(a) \quad \frac{\partial f_h}{\partial x}, \  \frac{\partial f_h}{\partial y} \ \text{exist and are continuous in $Q_h$.}
\\[8pt]
\qquad 
&(b) \quad \Gamma_{\un{Q}}[f_h] \ \leq \ \Gamma_{Q}[f]. 
\\[8pt]
\qquad 
&(c) \quad f_h \ra f \quad \ (h \ra 0) \ \text{uniformly in \un{Q}.}
\end{align*}

Granted these points, on the basis of the earlier considerations, from $(a)$, 

\[
\Gamma_{\un{Q}}[f_h] \ = \ L_{\un{Q}}[f_h],
\]
thus by $(b)$, 

\[
L_{\un{Q}}[f_h]  \ \leq \ \Gamma_{Q}[f] \ \leq L_Q[f]
\]

\qquad
$\implies$

\[
\limsup\limits_{h \ra 0} \ L_{\un{Q}}[f_h] \ \leq \ \Gamma_{Q}[f].
\]
But thanks to $(c)$,

\[
L_{\un{Q}}[f] \ \leq \ \liminf\limits_{h \ra 0} \ L_{\un{Q}}[f_h].
\]
And then
\begin{align*}
L_{\un{Q}}[f] \ 
&\leq \ \liminf\limits_{h \ra 0} \ L_{\un{Q}}[f_h].
\\[11pt]
&\leq \ 
\limsup\limits_{h \ra 0} \ L_{\un{Q}}[f_h]
\\[11pt]
&\leq  
\Gamma_{Q}[f]
\\[11pt]
&\leq \ L_Q[f].
\end{align*}
Suppose now that \un{$Q$} invades $Q:\un{Q} \uparrow Q$, hence
\[
L_{\un{Q}}[f] \ra L_Q[f]
\]

\qquad
$\implies$

\[
L_Q[f] \ \leq \ \Gamma_{Q}[f] \ \leq \ L_Q[f]
\]

\qquad
$\implies$

\[
\Gamma_{Q}[f] \ = \ L_Q[f].
\]

\chapter{
$\boldsymbol{\S}$\textbf{4}.\quad APPROXIMATION THEORY}
\setlength\parindent{2em}
\setcounter{theoremn}{0}
\renewcommand{\thepage}{\S4-\arabic{page}}

\qquad
To finish the proof that

\[
\Gamma_Q[f] = L_Q[f],
\]
we have yet to establish the validity of points $(a), (b), (c)$ as formulated near the end of the preceeding $\S$ and for this, it will be necessary to set up some machinery.
\\[-.25cm]

\begin{x}{\small\bf DEFINITION} \ 
Let $f:Q \ra \R$ be a continuous function and let $\ds 0 < h < \frac{1}{2}$ $-$then the function

\[
f_h(x,y) = \frac{1}{4h^2} \int_{-h}^h \int_{-h}^h f(x + \xi, y + \eta) 
\ \td\xi  \hsy \td\eta
\]
defined in the square

\[
Q_h: \ 
\begin{cases}
\ h \leq x \leq 1 - h
\\[4pt]
\ h \leq y \leq 1 - h
\end{cases}
\]
is called the \un{integral mean} of $f$.
\\[-.25cm]
\end{x}

\begin{x}{\small\bf LEMMA} \ 
$f_h:Q_h \ra \R$ is a continuous function.
\\[-.5cm]
\end{x}

\begin{x}{\small\bf LEMMA} \ 
$f_h \ra f$ $(h \ra 0)$ uniformly in $\un{Q} \subset Q_h$.
\\[-.25cm]
\end{x}

\begin{x}{\small\bf LEMMA} \ 
$\ds\frac{\partial f_h}{\partial x}$, 
$\ds\frac{\partial f_h}{\partial y}$ exist and are continuous functions on $Q_h:$

\[
\begin{cases}
\ \ds\frac{\partial f_h}{\partial x} 
\ = \  
\frac{1}{4h^2} \ 
\int_{-h}^h\  
[f(x + h, y + \eta) - f(x - h, y + \eta)]  
\ \td\eta 
\\[18pt]
\ \ds\frac{\partial f_h}{\partial y} 
\ = \  
\frac{1}{4h^2} \ 
\int_{-h}^h \ 
[f(x + \xi, y + h) - f(x+ \xi, y - h)] 
\ \td\xi
\end{cases}
.
\]
\\[-1cm]
\end{x}

\begin{x}{\small\bf \un{N.B.}}\ 
Accordingly points $(a)$ and $(c)$ are settled.
\\[-.5cm]
\end{x}

The validity of point $(b)$, i.e., the assertion that
\[
\Gamma_{\un{Q}}[f_h] \ \leq \ \Gamma_Q[f]
\]
is not so easy to prove.
\\[-.25cm]

Start by fixing an oriented rectangle $R \subset \un{Q}$: 

\[
\begin{cases}
\ a \leq x \leq b \qquad (a < b) 
\\[8pt]
\ c \leq y \leq d \qquad (c < d) 
\end{cases}
\abs{R} = (b - a)(d - c).
\]
Then
\[
\abs{f_h(x,d) - f_h(x,c)} \leq   \frac{1}{4h^2} \int_{-h}^h \int_{-h}^h \abs{f(x + \xi,d + \eta) - f(x + \xi,c + \eta)} \td\xi  \td\eta
\]

\noindent
$\implies$
\begin{align*}
G_\XX(f_h;R) \ 
&= \ \int_a^b \abs{f_h(x,d) - f_h(x,c)} 
\ \td x
\\[15pt]
&\leq \ 
\frac{1}{4h^2} \ 
\int_{-h}^h \ 
\int_{-h}^h  \ 
\ \td\xi  \hsy \td\eta 
\int_a^b \ 
\abs{f(x + \xi,d + \eta) - f(x + \xi,c + \eta)} 
\ \td x.
\end{align*}
Let $R_{\xi\eta}$ be the rectangle obtained by subjecting R to the translation
\[
\begin{cases}
\ \bar{x} = x + \xi 
\\[4pt]
\ \bar{y} = y + \eta, 
\end{cases}
\]
thus

\[
G_\X(f;R_{\xi\eta}) \ = \ 
\int_a^b \ 
\abs{f(x + \xi,d + \eta) - f(x + \xi,c + \eta)} 
\ \tdx
\]
and so

\[
G_\X(f_h;R) \ \leq \ 
\frac{1}{4h^2} \ 
\int_{-h}^h \ 
\int_{-h}^h \ 
G_\X(f;R_{\xi\eta}) 
\ \td\xi  \hsy \td\eta.
\]
Analogously

\[
G_\Y(f_h;R) \ \leq \ 
\frac{1}{4h^2} \ 
\int_{-h}^h \ 
\int_{-h}^h \ 
G_\Y(f;R_{\xi\eta}) 
\ \td\xi  \hsy \td\eta.
\]
Finally

\[
\abs{R} 
\ = \  
\abs{R_{\xi\eta}} 
\ = \ 
\frac{1}{4h^2} \ 
\int_{-h}^h \ 
\int_{-h}^h \ 
\abs{R_{\xi\eta}} 
\ \td\xi  \hsy \td\eta.
\]
To summarize:
\\[-.5cm]

\begin{x}{\small\bf LEMMA} \ 
\begin{align*}
\Gamma(f_h;R) \ 
&\leq \ 
[G_\X(f_h;R)^2 + G_\Y(f_h;R)^2 + \abs{R}^2]^{1/2}
\\[15pt]
&\leq \ 
\frac{1}{4h^2} \ 
\left[
\left(
\int_{-h}^h \ 
\int_{-h}^h  \ 
G_\X(f;R_{\xi\eta}) 
\ \td\xi  \hsy \td\eta
\right)^2 \right.
\\[15pt]
&\qquad + \left.
\left(
\int_{-h}^h \ 
\int_{-h}^h  \ 
G_\Y(f;R_{\xi\eta}) 
\ \td\xi  \hsy \td\eta
\right)^2 
+ 
\left(
\int_{-h}^h \ 
\int_{-h}^h  \ 
\abs{R_{\xi\eta}} 
\ \td\xi  \hsy \td\eta\right)^2
\right]^{1/2}.
\end{align*}
\\[-1cm]
\end{x}

\begin{x}{\small\bf RAPPEL} \ 
Under canonical assumptions,

\[
\left(\left(
\int_\X \ 
\phi_1\right)^2 
+ \cdots + 
\left(
\int_\X \ 
\phi_n\right)^2\right)^{1/2} 
\ \leq \ 
\int_\X \ 
\left(\phi_1^2 
+ \ldots + \phi_n^2\right)^{1/2}.
\]
\\[-1cm]

\noindent
Therefore
\begin{align*}
\Gamma(f_h;R) \ 
&\leq \ 
\frac{1}{4h^2} \ 
\int_{-h}^h \ 
\int_{-h}^h \ 
(G_\X(f;R_{\xi\eta})^2 + G_\Y(f;R_{\xi\eta})^2 + \abs{R_{\xi\eta}}^2)^{1/2}
\ \td\xi \hsy \td\eta
\\[15pt]
&= \ 
\frac{1}{4h^2} \int_{-h}^h \int_{-h}^h  \Gamma(f;R_{\xi\eta}) 
\ \td\xi \hsy \td\eta.
\\
\end{align*}
Suppose now that \un{$D$} is a subdivision of \un{$Q$} into nonoverlapping rectangles $R$ 
(lines parallel to the coordinate axes) $-$then
\begin{align*}
G(f;\un{D}) \ 
&= \ 
\sum \ 
\Gamma (f_h;R)
\\[15pt]
&\leq \ 
\frac{1}{4h^2} \ 
\int_{-h}^h \ 
\int_{-h}^h \ 
\sum \ 
\Gamma (f;R_{\xi\eta})
\ \td\xi \hsy \td\eta,
\end{align*}
the sum under $\ds\int_{-h}^h \int_{-h}^h$ being the sum of Ge\"ocze (for $f$) relative to the division 
$\un{D}_{\xi\eta}$ of $\un{Q}_{\xi\eta} \subset Q$ into rectangles $R_{\xi\eta}$, thus a fortiori, 

\[
\sum \ \Gamma (f;R_{\xi\eta}) \ \leq \ \Gamma_Q[f]
\]

\qquad
$\implies$
\begin{align*}
G(f_h;\un{D}) \ 
&\leq \  
\frac{1}{4h^2} \ 
\int_{-h}^h \ 
\int_{-h}^h  \ 
\Gamma_Q[f] 
\ \td\xi \hsy \td\eta
\\[15pt]
&= \ 
\frac{\Gamma_Q[f]}{4h^2} \ 
\int_{-h}^h \ 
\int_{-h}^h  \ 
\ \td\xi \hsy \td\eta
\\[15pt]
&= \ 
\Gamma_Q[f] 
\end{align*}

\qquad
$\implies$
\begin{align*}
\Gamma_{\un{Q}}[f_h]  \ 
&= \ 
\sup\limits_{\un{D}} \ G(f_h;\un{D})
\\[11pt]
&\leq \ 
\Gamma_Q[f],
\end{align*}
from which point $(b)$.
\\[-.25cm]
\end{x}

\begin{x}{\small\bf LEMMA} \ 

\[
L_{Q_h}[f_h] \ \leq \ L_Q[f]
\]
and
\[
L_Q[f] \ = \ \lim\limits_{h \ra 0} L_{Q_h}[f_h].
\]

Since
\[
L_{Q_h}[f_h] 
\ = \ 
\iint\limits_{Q_h} \ 
\left[1 + \left(\frac{\partial f_h}{\partial x}\right)^2 + \left(\frac{\partial f_h}{\partial y}\right)^2 
\hsy\right]^{1/2} 
\ \td x \td y,
\]
it follows that
\allowdisplaybreaks
\begin{align*}
L_Q[f] \ 
&= \ 
\lim\limits_{h \ra 0} \ 
\int_h^{1-h}
\int_h^{1-h} 
\left
[1 + 
\left(\frac{1}{4h^2} \ 
\int_{-h}^h \ 
(f(x+h,y+\eta) - f(x-h,y+\eta)) 
\ \td\eta
\right)^2 
\right.
\\[18pt]
& 
\qquad + \ \left.\left(\frac{1}{4h^2} \int_{-h}^h (f(x+\xi,y+h) - f(x+\xi, y-h)) 
\td\xi
\right)^2 \hsy \right]^{1/2} 
\ \td x \hsy \td y.
\end{align*}
\\[-1cm]
\end{x}

\[
*\ *\ *\ *\ *\ *\ *\ *\ *\ *\ *\ *\ 
\]
\\

What follows will not be needed in the sequel but it is of independent interest.
\\

\begin{x}{\small\bf DEFINITION} \ 
Let $f \in \Lp^1(Q)$ and let $\ds 0 < h < \frac{1}{2}$ $-$then the function
\[
f_h(x,y) 
\ = \ 
\frac{1}{4h^2} \ 
\int_{-h}^h  \ 
\int_{-h}^h \ 
f(x+\xi,y+\eta) 
\ \td\xi  \hsy \td\eta
\]
defined in the square
\[
\begin{cases}
\ h \leq x \leq 1 - h
\\[4pt]
\ h \leq y \leq 1 - h
\end{cases}
\]
is called the \un{integral mean} of $f$.
\\[-.25cm]
\end{x}

\begin{x}{\small\bf LEMMA} \ 
$f_h:Q_h \ra \R$ is a continuous function, hence

\[
\iint\limits_{Q_h} \abs{f_h} \ < +\infty \ \implies \ f_h \in \Lp^1(Q_h).
\]
\\[-1cm]
\end{x}

\begin{x}{\small\bf LEMMA} \ 
$\forall \ f \in \Lp^1(Q)$, 
\[
\norm{f_h}_{\Lp^1} \ \leq \ \norm{f}_{\Lp^1}.
\]

PROOF 
\allowdisplaybreaks
\begin{align*}
\iint\limits_{Q_h} \ \abs{f_h(x,y)} 
\ \td x \td y
&= \ 
\int_h^{1-h}\ 
\int_h^{1-h}\ 
 \abs{f_h(x,y)} 
\ \td x \td y
\\[18pt]
&\leq \ 
\frac{1}{4h^2}\  
\int_h^{1-h}\ 
\int_h^{1-h} \ 
\left\{\int_{-h}^h  \int_{-h}^h \abs{f(x+\xi, y+ \eta)} \ \td\xi \td\eta\right\}
\ \td x \td y
\\[18pt]
&\leq \ 
\frac{1}{4h^2}  \ 
\int_{-h}^h   \ 
\int_{-h}^h  \ 
\left\{
\int_h^{1-h} \ 
\int_h^{1-h}  \ 
\abs{f(x+\xi, y+ \eta)} 
\ \td x \td y
\right\} 
\ \td\xi \td\eta 
\\[18pt]
&\leq \ 
\frac{1}{4h^2}  \ 
\int_{-h}^h   \ 
\int_{-h}^{h}  \ 
\left\{
\int_{h+\xi}^{1-h+\xi}
\int_{h+\eta}^{1-h+\eta}
 \abs{f(x, y)} 
\ \td x \td y\right\} 
\td\xi \td\eta 
\\[18pt]
&\leq \ 
\frac{1}{4h^2} \  
\int_{-h}^h   \ 
\int_{-h}^{h}  \ 
\left\{
\int_{0}^{1} \ 
\int_{0}^{1}  \ 
\abs{f(x, y)} 
\ \td x \td y
\right\} 
\ \td\xi \td\eta 
\\[18pt]
&\leq \ 
\frac{1}{4h^2} \ (2h)(2h) \norm{f}_{\Lp^1}
\\[18pt]
&= \ 
\norm{f}_{\Lp^1}
\\[18pt]
&< \ 
+\infty.
\end{align*}
\\[-1cm]
\end{x}

\begin{x}{\small\bf REMARK} \ 
An analogous estimate obtains if $f \in \Lp^p(Q)$ $(1 < p < +\infty):$

\[
\norm{f_h}_{L^p} \ \leq \ \norm{f}_{\Lp^p}.
\]
\\[-1.25cm]
\end{x}

\begin{x}{\small\bf LEMMA} \ 
As $h \ra 0$, $f_h$ converges almost everywhere to $f$.
\\[-.5cm]
\end{x}

\begin{x}{\small\bf LEMMA} \ 

\[
\iint\limits_{Q_h} \abs{f_h - f} \ra 0 \qquad (h \ra 0).
\]

\begin{spacing}{1.75}
PROOF \ 
Given $\varepsilon > 0$, 
write $f = \phi + \psi$, where $\phi$ is continuous in $Q$, 
$\psi$ is integrable in $Q$, and $\ds \iint\limits_{Q} \abs{\psi} < \varepsilon$ $-$then
\end{spacing}
\begin{align*}
\iint\limits_{Q_h} \ \abs{f_h - f} \ 
&= \ 
\iint\limits_{Q_h} \ \abs{(\phi_h + \psi_h) - (\phi + \psi)}
\\[15pt]
&\leq \ 
\iint\limits_{Q_h} \ \abs{\phi_h - \phi} 
+ \iint\limits_{Q_h} \abs{\psi_h - \psi}
\\[15pt]
&\leq \
 \iint\limits_{Q_h} \ \abs{\phi_h - \phi} 
+ \iint\limits_{Q_h} \ \abs{\psi_h} 
+ \iint\limits_{Q_h} \ \abs{\psi}
\\[15pt]
&\leq \ 
\iint\limits_{Q_h} \ \abs{\phi_h - \phi} 
+ \iint\limits_{Q} \ \abs{\psi} 
+ \iint\limits_{Q} \ \abs{\psi}
\\[15pt]
&\leq \ 
\iint\limits_{Q_h} \ \abs{\phi_h - \phi} 
+ 2 \iint\limits_{Q} \ \abs{\psi}
\\[15pt]
&\leq \ 
\iint\limits_{Q_h} \ \abs{\phi_h - \phi} + 2\varepsilon.
\end{align*}
Since $\phi$ is continuous in $Q$, it follows that in $Q_h$,

\[
\phi_h \ra \restr{\phi}{{Q_h}} \qquad (h \ra 0)
\]
uniformly, hence

\[
\iint\limits_{Q_h} \ \abs{\phi_h - \phi} \ \ra 0 \qquad (h \ra 0).
\]
So for all sufficiently small $h$,

\[
\iint\limits_{Q_h} \ \abs{\phi_h - \phi} \ < \ \varepsilon
\]

\qquad
$\implies$

\[
\lim\limits_{h \ra 0} \ \iint\limits_{Q_h} \ \abs{f_h - f} \ < \ 3\varepsilon.
\]
\\[-.75cm]
\end{x}

\begin{x}{\small\bf REMARK} \ 
An analogous statement obtains if $f \in \Lp^p(Q)$ $(1 < p < +\infty):$
\[
\iint\limits_{Q_h} \abs{f_h - f}^p \ra 0 \qquad (h \ra 0).
\]
\\[-1cm]
\end{x}

\begin{x}{\small\bf LEMMA} \ 
If $f \in \Lp^p(Q)$ $(1 \leq p < +\infty)$, then
\[
\frac{\partial f_h}{\partial x} \quad \& \quad \frac{\partial f_h}{\partial y}
\]
belong to $\Lp^p(Q_h)$.
\\[-.5cm]

PROOF \ 
Take $p > 1$ and consider $\ds\frac{\partial f_h}{\partial x}$, thus

\[
\frac{\partial f_h}{\partial x} \ = \ \frac{1}{4h^2} \int_{y - h}^{y+h} f(x+h, \eta) - f(x-h, \eta)  \td\eta
\]
almost everywhere in $Q_h$, the claim being that the functions

\[
\begin{cases}
\ \ds
\int_{y - h}^{y+h} \ f(x+h, \eta)  \ \td\eta 
\\[26pt]
\ \ds
\int_{y - h}^{y+h} \ f(x-h, \eta)  \ \td\eta 
\end{cases}
\]
are in $\Lp^p(Q_h)$.  To discuss the first of these, write

\[
\int_{y - h}^{y+h} \ f(x+h, \eta)  \ \td\eta 
\ = \  
\int_{-h}^h f(x+h,y+\eta)  \ \td\eta.
\]

\noindent
Then

\[
\abs{\int_{-h}^h f(x+h,y+\eta)  
\ \td\eta}^p 
\ \leq \ 
(2h)^{p-1} \ 
\int_{-h}^h \ 
\abs{f(x+h,y+\eta)}^p  
\ \td\eta.
\]
Since $f \in \Lp^p(Q)$, $\abs{f(x+h,y+\eta)}^p$ is integrable in

\[
h \leq x \leq 1 - h , \quad h \leq y \leq 1 - h, \quad -h \leq \eta \leq h.
\]
Therefore

\[
\int_{-h}^h \ 
\abs{f(x+h,y+\eta)}^p 
\ \td\eta
\]
is integrable in $Q_h$, hence

\[
\int_{-h}^h f(x+h,y+\eta)  \td\eta
\]
is in $\Lp^p(Q_h)$.
\end{x}

\chapter{
$\boldsymbol{\S}$\textbf{5}.\quad TONELLI'S CHARACTERIZATION}
\setlength\parindent{2em}
\setcounter{theoremn}{0}
\renewcommand{\thepage}{\S5-\arabic{page}}

\vspace{-.5cm}
\qquad
Let $f:Q \ra \R$ be a continuous function.
\\[-.5cm]

\begin{x}{\small\bf DEFINITION} \ 
\[
\begin{cases}
\ \tV_x(f;y) \ = \ T_{f(-,y)} [0,1] \qquad (0 \leq y \leq 1) 
\\[4pt]
\ \tV_y(f;x) \ = \ T_{f(x,-)} [0,1] \qquad (0 \leq x \leq 1)
\end{cases}
.
\]
\\[-1cm]
\end{x}

\begin{x}{\small\bf LEMMA} \ 
\[
\begin{cases}
\ \tV_x(f;-) \text{ is a lower semicontinuous function of $y \in [0,1]$} 
\\[4pt]
\ \tV_y(f;-) \text{ is a lower semicontinuous function of $x \in [0,1]$} 
\end{cases}
.
\]

PROOF \ 
Consider the first assertion and suppose that $y_n \ra y$ $-$then
\[
f(x,y_n) \ra f(x,y)
\]

\qquad
$\implies$
\[
T_{f(-,y)}[0,1] \ \leq \ \liminf\limits_{n \ra \infty} \ T_{f(-,y_n)}[0,1].
\]
I.e.:
\[
\tV_x(f;y) \ \leq \ \liminf\limits_{n \ra \infty} \ \tV_x(f;y_n).
\]
\\[-1cm]
\end{x}

\begin{x}{\small\bf SCHOLIUM} \ 
$V_x(f;-)$ and $V_y(f;-)$ are Lebesgue measurable.
\\[-.5cm]
\end{x}

\begin{x}{\small\bf DEFINITION} \ 
(BVT) $f$ is said to be of \un{bounded variation in the sense} \un{of Tonelli} if
\[
\begin{cases}
\ \ds
\int\limits_0^1 \ 
\tV_x(f;y) \tdy \ < \  +\infty  
\\[18pt]
\ \ds
\int\limits_0^1 \ 
\tV_y(f;x) \tdx \ < \  +\infty
\end{cases}
.
\]
\\[-1cm]
\end{x}

\begin{x}{\small\bf NOTATION} \ 
\[
\tV_T(f) 
\ = \ 
\int\limits_0^1 \ \tV_x(f;y) \ \tdy \hsx + \hsx \int\limits_0^1\  \tV_y (f;x) \ \tdx.
\]
\\[-1cm]
\end{x}

\begin{x}{\small\bf \un{N.B.}} \ 
Accordingly, if $\tV_T(f) < +\infty$, then
\[
e_Y \ = \ \{y \in [0,1]: \tV_x(f;y) = +\infty\}
\]
is of Lebesgue measure zero and
\[
e_X \ = \ \{x \in [0,1]: \tV_y(f;x) = +\infty\}
\]
is of Lebesgue measure zero.
\\[-.5cm]
\end{x}

\begin{x}{\small\bf LEMMA} \ 
Suppose that $\tV_T(f) < +\infty$ $-$then $\restr{f}{Q^\circ} \in \BV (Q^\circ)$ and 
\[
\begin{cases}
\ \ds
f_x \ = \ \frac{\partial f}{\partial x} \text{ exists almost everywhere in $Q$}  
\\[11pt]
\ \ds
f_y \ = \ \frac{\partial f}{\partial y} \text{ exists almost everywhere in $Q$}  
\end{cases}
.
\]
\\[-1cm]
\end{x}

\begin{x}{\small\bf LEMMA} \ 
Suppose that $\tV_T(f) < +\infty$ $-$then
\[
\begin{cases}
\ \ds
\iint\limits_Q \ 
\abs{f_x(x,y)} \ 
\tdx \tdy \ \leq \ \int\limits_0^1 \tV_x(f;y) \ 
\tdy < +\infty  
\\[18pt]
\ \ds
\iint\limits_Q \ 
\abs{f_y(x,y)} \ 
\tdx \tdy \ \leq \ \int\limits_0^1 \tV_y(f;x) \ 
\tdx < +\infty  
\end{cases}
\]
\\[-1cm]

\qquad
$\implies$
\[
\begin{cases}
\ f_x  
\\[4pt]
\ f_y  
\end{cases}
\in L^1(Q)
\]

\qquad
$\implies$
\[
\left[1 + f_x^2 + f_y^2\right]^{1/2} \in L^1(Q).
\]
\\[-1cm]
\end{x}

\begin{x}{\small\bf THEOREM} \ 
$L_Q[f]$ is finite iff $f$ is of bounded variation in the sense of Tonelli.
\\[-.25cm]

Assume to begin with that $L_Q[f]$ is finite.  
Let $D$ be the subdivision of $Q$ specified by
\[
\begin{cases}
\ x_0 = 0 < x_1 < \ldots < x_j < \ldots < x_m = 1
\\[4pt]
\ y_0 = 0 < y_1 < \ldots < y_k < \ldots < y_n = 1
\end{cases}
\]
and introduce

\[
\begin{cases}
\ \ds 
\textnormal{v}_x(f;y;D) 
\ = \ 
\sum\limits_{j=0}^{m-1}\ 
\abs{f(x_{j+1},y) - f(x_j,y)} \qquad (0 \leq y \leq 1)
\\[18pt]
\ \ds
\textnormal{v}_y(f;x;D) 
\ = \ 
\sum\limits_{k=0}^{n-1} \ 
\abs{f(x,y_{k+1} - f(x,y_k)} \qquad (0 \leq x \leq 1)
\end{cases}
.
\]
Then

\[
\begin{cases}
\ \ds
\int\limits_0^1 \ 
\textnormal{v}_x(f;y;D) \ \tdy 
\ = \ 
\sum \ G_Y(f;R)
\\[18pt]
\ \ds
\int\limits_0^1 \ 
\textnormal{v}_y(f;x;D) \ \tdx 
\ = \ 
\sum \ G_X(f;R)
\end{cases}
,
\]
\\[-.75cm]

\noindent
the summations being over the rectangles $R$ in $D$. 
Next
\[
\begin{cases}
\ \ds
\sum \ G_Y(f;R)
\\[15pt]
\ \ds
\sum \ G_X(f;R)
\end{cases}
\ \ \leq \ G(f;D) \ \leq\  L_Q[f].
\]
Therefore
\[
\begin{cases}
\ \ds
\int\limits_0^1 \ 
\textnormal{v}_x(f;y;D) \ \tdy
\\[26pt]
\ \ds
\int\limits_0^1 \ 
\textnormal{v}_y(f;x;D) \ \tdx
\end{cases}
\ \leq\  L_Q[f] < +\infty.
\]
\\[-.75cm]

\noindent
From the definitions,

\[
\begin{cases}
\ 0 \ \leq \ \textnormal{v}_x(f;y;D) \ \leq \ \tV_x(f;y)
\\[11pt]
\ 0 \ \leq \  \textnormal{v}_y(f;x;D) \ \leq \ \tV_y(f;x)
\end{cases}
.
\]
So, upon sending the maximum diameters of the rectangles $R$ in $D$ to zero sequentially, we conclude that

\[
\begin{cases}
\ \lim \textnormal{v}_x(f;y;D) \ \leq \ \tV_x(f;y)
\\[11pt]
\ \lim \textnormal{v}_y(f;x;D) \ \leq \ \tV_y(f;x)
\end{cases}
\]

\qquad
$\implies$

\[
\begin{cases}
\ \ds
\int\limits_0^1 \ 
\tV_x(f;y) 
\ \tdy 
\ = \ 
\int\limits_0^1  \ 
\lim \ 
\textnormal{v}_x(f;y; D) \ \tdy 
\\[18pt]
\ \ds
\int\limits_0^1 \ 
\tV_y(f;x) 
\ \tdx 
\ = \ 
\int\limits_0^1  \ 
\lim \ 
\textnormal{v}_y(f;x; D) \ \tdx 
\end{cases}
\]
or still, 

\[
\begin{cases}
\ \ds
\leq \ 
\liminf \int\limits_0^1  \ 
\textnormal{v}_x(f;y;D) \ \tdy 
\\[18pt]
\ \ds 
\leq \ 
\liminf \int\limits_0^1  \ 
\textnormal{v}_y(f;x;D) \ \tdx 
\end{cases}
(\text{Fatou}) \quad \leq L_Q[f] \ < \ +\infty.
\]

\noindent
Consequently, under the supposition that $L_Q[f]$ is finite, 
it follows that $f$ is of bounded variation in the sense of Tonelli.
\\[-.25cm]

To reverse this, note first that for any $D$,

\[
\begin{cases}
\ \textnormal{v}_x(f;y;D) \ \leq \ \tV_x(f;y)
\\[11pt]
\ \textnormal{v}_y(f;x;D) \ \leq \ \tV_y(f;x)
\end{cases}
\]
\qquad
$\implies$
\vspace{0.25cm}
\[
\begin{cases}
\ \ds 
\sum \ 
G_Y(f;R) 
\ \leq \ 
\int\limits_0^1 \ 
\tV_x(f;y)\ \tdy
\\[15pt]
\ \ds
\sum G_X(f;R) 
\ \leq \ 
\int\limits_0^1 \ 
\tV_y(f;x)\ \tdx
\end{cases}
.
\]
\\[-.75cm]

\noindent
And
\[
\Gamma(f;R) \ \leq \ G_X(f;R) + G_Y(f;R) + \abs{R}
\]
\qquad
$\implies$
\begin{align*}
G(f;D) \ 
&= \ 
\sum \ 
\Gamma(f;R)
\\[18pt]
&\leq \ 
\sum \ 
G_Y(f;R) + \sum \ G_X(f;R) + \sum \ \abs{R}
\\[18pt]
&\leq \ 
\int\limits_0^1 \ 
\tV_x(f;y) \ \tdy + \int\limits_0^1 \ \tV_y(f;x) 
\ \tdx + 1
\\[18pt]
&= \ 
\tV_T(f) + 1.
\end{align*}
However
\[
\Gamma_Q[f] 
\ = \  
\sup\limits_{D} \ G(f;D).
\]
Therefore
\[
\Gamma_Q[f] < +\infty
\]
\qquad
$\implies$
\[
L_Q[f] 
\ < \  
+\infty.
\]
\\[-1cm]
\end{x}

\begin{x}{\small\bf REMARK} \ 
Individually
\[
\int\limits_0^1 \ \tV_x(f;y) \ \tdy,
\quad  
\int\limits_0^1 \ \tV_y(f;x) \ \tdx,\  1
\]
are all $\leq L_Q[f]$.
\end{x}

\chapter{
$\boldsymbol{\S}$\textbf{6}.\quad TONELLI'S ESTIMATE}
\setlength\parindent{2em}
\setcounter{theoremn}{0}
\renewcommand{\thepage}{\S6-\arabic{page}}

\qquad
Let $f:Q \rightarrow \R$ be a continuous function.
\\[-0.25cm]

\begin{x}{\small\bf THEOREM} \ 
Suppose that $L_Q[f]$ is finite $-$then
\[
L_Q[f] \ \geq \ \iint\limits_Q \left[1 + f_x^2 + f_y^2\right]^{1/2} \ \tdx \hsy \tdy.
\]
\\[-1.5cm]
\end{x}

Let $D = \{R_1, R_2, \ldots, R_n\}$ be a subdivision of $Q$, where
\[
R_k \ = \ [a_k,b_k] \times [c_k,d_k] \qquad (k = 1, 2, \ldots, n).
\]
\\[-1.5cm]

\begin{x}{\small\bf LEMMA} \ 
Given $\varepsilon > 0$, there is a $D$ such that
\begin{align*}
\bigg|
\sum\limits_{k=1}^n \ 
\bigg[
\bigg(
\iint\limits_{R_k} \ 
f_x 
\ \tdx  \hsy \tdy
\bigg)^2 
&+ 
\bigg(
\iint\limits_{R_k} \ 
f_y 
\ \tdx \hsy \tdy
\bigg)^2 
+ \abs{R_k}^2\bigg]^{1/2} 
\\[18pt]
&-  
\iint\limits_Q \ 
\bigg[
1 + f_x^2 + f_y^2
\bigg]^{1/2} 
\ \tdx \hsy \tdy \hsx
\bigg|
\ < \  
\varepsilon.
\end{align*}

[Recall that
\[
\begin{cases}
\ f_x
\\[4pt]
\ f_y
\end{cases}
\ \in L^1(Q)
\]
\\[-.75cm]

\noindent
and use the Vitali covering lemma.]
\\[-0.25cm]

Proceeding
\[
\bigg|
\sum\limits_{k=1}^n\  [\ldots]^{1/2} - \iint\limits_Q \ldots
\bigg|
\ < \  
\varepsilon
\]

\qquad
$\implies$

\[
\bigg|
\iint\limits_Q \ -
 \sum\limits_{k=1}^n \ 
 [\ldots]^{1/2}
 \bigg|
 \ < \   
 \varepsilon
\]

\qquad
$\implies$

\[
\iint\limits_Q 
\ - \  
\sum\limits_{k=1}^n \ [\ldots]^{1/2} \ < \   \varepsilon
\]

\qquad
$\implies$

\[
\sum\limits_{k=1}^n \ [\ldots]^{1/2} 
\ - \ 
\iint\limits_Q \ \cdots 
\ > \   -\varepsilon
\]

\qquad
$\implies$

\[
\sum\limits_{k=1}^n [\ldots]^{1/2} 
\ > \ 
\iint\limits_Q \ \cdots - \varepsilon.
\]
And
\[
\Gamma_Q[f] 
\ \geq \ 
\sum\limits_{k=1}^n \ [\ldots]^{1/2} 
\ > \ 
\iint\limits_Q \ \cdots - \varepsilon.
\]
But
\[
\Gamma_Q[f] 
\ = 
\ L_Q[f].
\]
\end{x}

\chapter{
$\boldsymbol{\S}$\textbf{7}.\quad THE ROLE OF ABSOLUTE CONTINUITY}
\setlength\parindent{2em}
\setcounter{theoremn}{0}
\renewcommand{\thepage}{\S7-\arabic{page}}

\vspace{-0.5cm}
\qquad 
Let $f: Q \ra \R$ be a continuous function.
\\[-.5cm]

\begin{x}{\small\bf DEFINITION} \ (ACT) \ 
$f$ is said to be 
\un{absolutely continuous in the sense} \un{of Tonelli} 
if it is of bounded variation in the sense of Tonelli and if 
\[
\begin{cases}
\ \text{For almost every $y \in [0,1]$, the function 
$x \ra f(x,y)$ is absolutely continuous}
\\[4pt]
\ \text{For almost every $x \in [0,1]$, the function 
$y\ra f(x,y)$ is absolutely continuous}
\end{cases}
.
\]
\end{x}

\begin{x}{\small\bf REMARK} \ 
Since $f$ is BVT, the ordinary partial derivatives
\[
\frac{\partial f}{\partial x} 
\ \ \& \ \
\frac{\partial f}{\partial y} 
\]
belong to $\Lm^1(Q)$.  
So, thanks to ACL, 
\[
f \in W^{1, 1} (Q^\circ).
\]
\end{x}

\begin{x}{\small\bf NOTATION} \ 
Put 
\[Q^{(h, k)} 
\ = \ 
[0, 1 - h] \hsx \times \ \hsx [0, 1 - k],
\]
where
\[
\begin{cases}
\ 0 < h < 1
\\[4pt]
\ 0 < k < 1
\end{cases}
.
\]
\\[-.9cm]
\end{x}

\begin{x}{\small\bf PICTURE} \ 
\\[-.9cm]
\[
\begin{tikzpicture}[scale=1.7]

      \draw[-] (-1,0) -- (2,0) node[right] {$$};
      \node[label={{}}] at (0,-.25) {$0$};
      \node[label={{}}] at (1,-.25) {$1-h$};
      \node[label={{}}] at (2,-.25) {$1$};
      
	\draw[-]  (0,0) -- (0,2) node[right] {$$};
      \node[label={{}}] at (-.4,1) {$1 - k$};
      \node[label={{}}] at (-.13,2) {$1$};
      
     \draw[-]  (0,1) -- (1,1) node[right] {$$};
     \draw[-]  (1,0) -- (1,1) node[right] {$$};
     
     \node[label={{}}] at (0.5,0.5) {$Q^{(h, k)}$};
     \draw[-]  (0,2) -- (2,2) node[right] {$$};
     \draw[-]  (2,0) -- (2,2) node[right] {$$};
      
      

\end{tikzpicture}
\]
\end{x}

\begin{x}{\small\bf NOTATION} \ 
Given an ACT function $f$, put
\[
f^{(h, k)} (x, y) 
\ = \ 
\frac{1}{h \hsy k} \ 
\int\limits_x^{x + h} \ 
\int\limits_y^{y + k} \ 
f(\xi, \eta) \ 
\td \xi \hsy \td \eta.
\]
\end{x}

\begin{x}{\small\bf LEMMA} \ 
\[
\int\limits_0^{1 - h} \ 
\int\limits_0^{1 - k} \ 
\abs{f^{(h, k)} (x, y) }
\td x \hsy \td y \ 
\ \leq \ 
\int\limits_0^1 \ 
\int\limits_0^1 \ 
\abs{f(x,y)}\ 
\td x \hsy \td y.
\]
\\[-.75cm]
\end{x}

\begin{x}{\small\bf LEMMA} \ 
\[
\begin{cases}
\ \ds \frac{\partial f^{(h, k)}}{\partial x} 
\ = \ 
\frac{1}{h k} \ 
\int\limits_x^{x + h} \ 
\int\limits_y^{y + k} \ 
\frac{\partial f}{\partial \xi} 
\ \td \xi \hsy \td \eta
\\[26pt]
\ \ds \frac{\partial f^{(h, k)}}{\partial y} 
\ = \ 
\frac{1}{h k} \ 
\int\limits_x^{x + h} \ 
\int\limits_y^{y + k} \ 
\frac{\partial f}{\partial \eta} 
\ \td \xi \hsy \td \eta
\end{cases}
.
\]
\\[-1cm]

[Note: \ 
It therefore follows from these relations that $f^{(h, k)}$ is a $C^\prime$ function.]
\end{x}

Therefore
\allowdisplaybreaks
\begin{align*}
\int\limits_0^{1 - h} \ 
\int\limits_0^{1 - k} \ 
&
\sqrt{1 + \big[f_x^{(h, k)}\big]^2 + \big[f_y^{(h, k)}\big]^2} \ 
\td x \hsy  \td y \ 
\\[15pt]
&\hspace{.5cm}
=\ 
\int\limits_0^{1 - h} \ 
\int\limits_0^{1 - k} \ 
\bigg\{
\sqrt
{
\bigg[
\frac{1}{h k} \ 
\int\limits_0^h \ 
\int\limits_0^k \
\td \xi \hsy \td \eta 
\bigg]^2
}
\\[15pt]
&\hspace{2cm}
\ov{
+ 
\bigg[
\frac{1}{h k} \ 
\int\limits_0^h \ 
\int\limits_0^k \
f_\xi (x + \xi, y + \eta)
\ \td \xi \hsy \td \eta 
\bigg]^2
}
\\[15pt]
&\hspace{2cm}
\ov{
+ 
\bigg[
\frac{1}{h k} \ 
\int\limits_0^h \ 
\int\limits_0^k \
f_\eta (x + \xi, y + \eta)
\ \td \xi \hsy \td \eta 
\bigg]^2
\bigg\}
}
\ \td x \hsy \td y  
\\[15pt]
&\hspace{.5cm}
\leq\ 
\int\limits_0^{1 - h} \ 
\int\limits_0^{1 - k} \ 
\bigg[
\frac{1}{h k} \ 
\int\limits_0^h \ 
\int\limits_0^k \
\bigg\{
\sqrt{1 + \big[f_\xi (x + \xi, y + \eta)\big]^2
}
\\[15pt]
&\hspace{3cm}
\ov{+  \quad 
\big[f_\eta (x + \xi, y + \eta)\big]^2
\bigg\}
\ \td \xi \hsy \td \eta 
\bigg]} 
\ \td x \hsy \td y 
\\[15pt]
&\hspace{.5cm}
=\ 
\frac{1}{h k} \ 
\int\limits_0^h \ 
\int\limits_0^k \
\bigg[
\int\limits_\xi^{1 - h + \xi} \ 
\int\limits_\eta^{1 - k + \eta} \ 
\sqrt{1 + f_x^2 + f_y^2} 
\ \td x \hsy \td y 
\bigg] 
\ \td \xi \hsy \td \eta 
\\[15pt]
&\hspace{.5cm}
\leq\ 
\frac{1}{h k} \ 
\int\limits_0^h \ 
\int\limits_0^k \
\bigg[
\int\limits_0^1 \ 
\int\limits_0^1 \ 
\sqrt{1 + f_x^2 + f_y^2} 
\ \td x \hsy \td y 
\bigg]  
\ \td \xi \hsy \td \eta 
\\[15pt]
&\hspace{.5cm}
=\ 
\frac{1}{h k} \ 
\frac{h k}{1} \ 
\int\limits_0^1 \ 
\int\limits_0^1 \ 
\sqrt{1 + f_x^2 + f_y^2} 
\ \td x \hsy \td y 
\\[15pt]
&\hspace{.5cm}
=\ 
\iint\limits_Q \ 
\big[1 + f_x^2 + f_y^2 \big]^{1/2} \ 
\td x \td y
\\[15pt]
&\hspace{.5cm}
\leq\ 
L_Q [f].
\end{align*}
\\[-.75cm]

\begin{x}{\small\bf RAPPEL} \ 
During the course of establishing that

\[
\Gamma_Q [f]
\ = \ 
L_Q [f],
\]
it was shown that if $f$ was $C^\prime$, then 
\[
L_Q [f] 
\ = \ 
\iint\limits_Q \ 
\bigg[
1 + \bigg(\frac{\partial f}{\partial x}\bigg)^2 + \bigg(\frac{\partial f}{\partial y}\bigg)^2
\hsx \bigg]^{1/2} 
\ \td x \hsy \td y.
\]

So, upon applying this to $f^{(h, k)}$, the upshot is that
\[
L_{Q^{(h,k)}} \big[f^{(h,k)}\big]
\ = \ 
\int\limits_0^{1 - h} \ 
\int\limits_0^{1 - k} \ 
\sqrt{
\big[1 + [f_x^{(h,k)}\big]^2 + \big[f_y^{(h,k)}\big]^2 
}
\ \td x \hsy \td y.
\]
\end{x}

\begin{x}{\small\bf SCHOLIUM} \ 
If $f$ is absolutely continuous in the sense of Tonelli, then

\[
L_Q [f] 
\ = \ 
\iint\limits_Q \ 
\big[1 + f_x^2 + f_y^2\big]^{1/2} 
\ \td x \hsy \td y.
\]

[In fact, 
\allowdisplaybreaks
\begin{align*}
L_Q [f] \ 
&\leq \ 
\liminf\limits_{\substack{h \ra 0\\k \ra 0}} \ 
L_{Q^{(h, k)}} \hsx [f^{(h, k)}]
\\[15pt]
&\leq 
\limsup\limits_{\substack{h \ra 0\\k \ra 0}} \ 
L_{Q^{(h, k)}} \hsx \big[f^{(h, k)}\big]
\\[15pt]
&\leq \ 
\iint\limits_Q \ 
\big[1 + f_x^2 + f_y^2 \big]^{1/2} 
\ \td x \hsy \td y
\\[15pt]
&\leq \
L_Q [f] \hsy.]
\end{align*}
\\[-1.25cm]
\end{x}

\begin{x}{\small\bf EXAMPLE} \ 
Suppose that $f: \R^2 \ra \R$ is a $C^\prime$ function.  
Put

\[
\Gr_f(Q)
\ = \ 
\{(x, y), f(x, y) : (x, y) \in Q\}.
\]
Then
\[
\sH^2 (\Gr_f(Q)) 
\ = \ 
\iint\limits_Q \ 
\big[1 + f_x^2 + f_y^2 \big]^{1/2} 
\ \td x \hsy \td y.
\]
Consequently

\[
\sH^2 (\Gr_f(Q)) 
\ = \ 
L_Q [f].
\]
\end{x}

Matters can be reversed, namely: 
\\[-.25cm]

\begin{x}{\small\bf SCHOLIUM} \ 
If $f$ is of bounded variation in the sense of Tonelli and if 
\[
L_Q [f]
\ = \ 
\iint\limits_Q \ 
\big[1 + f_x^2 + f_y^2\big]^{1/2}  
\ \td x \hsy \td y,
\]
then $f$ is absolutely continuous in the sense of Tonelli.
\end{x}

We shall sketch the proof.
\\[-.25cm]

\begin{x}{\small\bf LEMMA} \ 
For every oriented rectangle $R \subset Q$, 
\[
L_R [f]
\ = \ 
\iint\limits_R \ 
\big[1 + f_x^2 + f_y^2 \big]^{1/2} 
\ \td x \hsy \td y.
\]
\end{x}

Explicate $R \subset Q$:

\[
\begin{cases}
\ a \leq x \leq b \quad (a < b)
\\[8pt]
\ c \leq y \leq d \quad (c <d)
\end{cases}
, \ 
\abs{R} = (b - a) \hsx (d - c)
\]
and introduce

\[
\begin{cases}
\ds
\ W_x (f; R) 
\ = \ 
\int\limits_c^d \ 
\tV_x (f;y) 
\ \td y
\\[18pt]
\ds
\ W_y (f; R)
\ = \ 
\int\limits_a^b \ 
\tV_y (f;x) 
\ \td x
\end{cases}
.
\]
\\[-.75cm]

\begin{x}{\small\bf LEMMA} \ 
For every oriented rectangle $R \subset Q$, 
\[
\begin{cases}
\ \ds W_x (f; R) 
\ \leq \ 
L_R [f]
\\[8pt]
\ \ds  W_y (f; R)
\ \leq \ 
L_R [f]
\end{cases}
.
\]
\end{x}

Therefore

\[
\begin{cases}
\ \ds W_x (f; R) 
\ \leq \ 
\iint\limits_R \ 
\big[1 + f_x^2 + f_y^2 \big]^{1/2}  
\ \td x \hsy \td y
\\[26pt]
\ \ds  W_y (f; R)
\ \leq \ 
\iint\limits_R \ 
\big[1 + f_x^2 + f_y^2\big]^{1/2}  
\ \td x \hsy \td y
\end{cases}
.
\]
\\[-.5cm]

Denoting by $\sR$ the set of oriented rectangles in $Q$, a 
\un{rectangle function} is a function 
$\phi:\sR \ra \R$.  
So, e.g., the assignments
\[
\begin{cases}
\ R \ra W_x (f; R)
\\[4pt]
\ R \ra W_y (f; R)
\end{cases}
(R \in \sR)
\]
are rectangle functions.
\\

\begin{x}{\small\bf DEFINITION} \ 
A rectangle function $R \ra \phi(R)$ is said to be 
\un{absolutely} \un{continuous} 
if for every $\varepsilon > 0$ there exists $\delta > 0$ such that 
\[
\abs{\phi(R_1)} + \cdots + \abs{\phi(R_n)} 
\ < \  
\varepsilon
\]
for every finite system of oriented rectangles 
$R_1, \ldots, R_n$ 
which satisfy the conditions
\[
R_i^\circ \cap R_j^\circ = \emptyset 
\ (i \neq j) 
\quad \text{and} \quad 
\abs{R_1} + \cdots + \abs{R_n} 
\ < \   
\delta.
\]
\\[-1.25cm]
\end{x}

\begin{x}{\small\bf CRITERION} \ 
If $\Phi \in \Lm^1 (Q)$ and if 
\[
\phi (R) 
\ = \ 
\iint\limits_R \ 
\abs{\Phi} 
\ \td x \hsy \td y
\quad (R \in \sR),
\]
then $\phi$ is absolutely continuous.
\\[-.5cm]
\end{x}

\begin{x}{\small\bf APPLICATION} \ 
The rectangle functions
\[
\begin{cases}
\ R \ra W_x (f; R)
\\[4pt]
\ R \ra W_y (f; R)
\end{cases}
(R \in \sR)
\]
are absolutely continuous.
\\[-.5cm]

[Note: \ 
Bear in mind that
\[
\big[1 + f_x^2 + f_y^2 \big]^{1/2} \in \Lm^1 (Q).]
\]
\\[-.75cm]

Recall that the contention is that $f$ is absolutely continuous in the sense of Tonelli, i.e., 
\[
\begin{cases}
\ \text{For almost every $y \in [0,1]$, the function 
$x \ra f(x,y)$ is absolutely continuous}
\\[4pt]
\ \text{For almost every $x \in [0,1]$, the function 
$y \ra f(x,y)$ is absolutely continuous}
\end{cases}
.
\]
\\[-1cm]
\end{x}

Consider the first of these assertions.  
Using the absolute continuity of 
$W_x (f; R)$
to eliminate a potential singular term, we have

\[
W_x (f; Q)
\ = \ 
\iint\limits_Q \ 
\abs{f_x(x, y)}
\ \td x \hsy \td y.
\]
On the other hand, by definition, 

\[
W_x (f; Q)
\ = \
\int\limits_0^1 \ 
\tV_x (f; y) 
\ \td y.
\]
Therefore

\[
\int\limits_0^1 \ 
[\tV_x (f; y) - 
\int\limits_0^1 \ \abs{f_x(x, y)}\ 
\td x \hsy ] \hsy 
\td y 
\ = \ 
0.
\]
But

\[
\tV_x (f; y) 
\ \geq \ 
\int\limits_0^1 \ 
\abs{f_x(x, y)}\ 
\td x
\]
for almost every $y \in [0,1]$.  
Therefore
\[
\tV_x (f; y) 
\ = \ 
\int\limits_0^1 \ 
\abs{f_x(x, y)}
\ \td x
\]
for those $y \notin E$, where $E$ is a certain subset of $[0,1]$ of Lebesgue measure 0.  
And this implies that $f(x,y)$ is absolutely continuous as a function of $x$ for $y \notin E$.
\\

\begin{x}{\small\bf \un{N.B.}}  \ 
In general, if $f$ is of bounded variation in the sense of Tonelli, then
\allowdisplaybreaks
\begin{align*}
W_x (f; Q) \ 
&=\ 
\int\limits_0^1 \ \tV_x (f; y) 
\ \td y
\\[15pt]
&\geq \ 
\int\limits_0^1 \
\bigg[
\int\limits_0^1 \
\abs{f_x(x, y)} 
\ \td x
\bigg]
\td y
\\[15pt]
&= \ 
\iint\limits_Q \ 
\abs{f_x(x, y)}
\ \td x \hsy \td y,
\end{align*}
the inequality becoming an equality in the presence of the absolute continuity of 
$R \ra W_x (f; R)$.
\end{x}

\chapter{
$\boldsymbol{\S}$\textbf{8}.\quad STEINER'S INEQUALITY}
\setlength\parindent{2em}
\setcounter{theoremn}{0}
\renewcommand{\thepage}{\S8-\arabic{page}}

\qquad 
Suppose that

\[
\begin{cases}
\ f_1 : Q \ra \R
\\[4pt]
\ f_2 : Q \ra \R
\end{cases}
\]
are continuous functions.
\\

\begin{x}{\small\bf THEOREM} \ 

\[
L_Q [(f_1 + f_2)/2]
\ \leq \ 
\frac{L_Q [f_1] + L_Q [f_2]}{2}.
\]

PROOF \ 
The assertion is trivial if 

\[
L_Q [f_1] = +\infty 
\quad \text{or} \quad
L_Q [f_2] = +\infty, 
\]
so it can be assumed that both are finite.  
Accordingly, given a subdivision $D$ of $Q$, form the sums of the Ge\"ocze per $f_1$, $f_2$, and 
$(f_1 + f_2)/2$, hence

\[
G((f_1 + f_2)/2; D) 
\ \leq \ 
\frac{G (f_1; D) + G (f_2; D)}{2}
\]

\qquad 
$\implies$ 

\[
G((f_1 + f_2)/2; D) 
\ \leq \ 
\frac{L_Q [f_1]  + L_Q [f_2] }{2}
\]

\qquad 
$\implies$ 

\[
L_Q [(f_1 + f_2)/2]
\ \leq \ 
\frac{L_Q [f_1]  + L_Q [f_2] }{2}.
\]
\\[-1cm]
\end{x}

\begin{x}{\small\bf RAPPEL} \ 
If $f: Q \ra \R$ is continuous, then $L_Q [f]$ is finite iff $f$ is of bounded variation in the sense of Tonelli, 
there being the estimate

\[
L_Q [f]
\ \geq \ 
\iint\limits_Q \ 
\big[1 + f_x^2 + f_y^2\big]^{1/2} \ 
\td x \hsy \td y, 
\]
the inequality becoming an equality iff $f$ is absolutely continuous in the sense of Tonelli.
\\

Suppose that

\[
\begin{cases}
\ f_1 : Q \ra \R
\\[4pt]
\ f_2 : Q \ra \R
\end{cases}
\]
are absolutely continuous in the sense of Tonelli $-$then 
the same is true of $(f_1 + f_2)/2$ and Steiner's inequality is the relation

\begin{align*}
&\iint\limits_Q \ 
\bigg\{\frac
{
\big[
1 + f_{1 x}^2 + f_{1 y}^2
\big]^{1/2}
+
\big[
1 + f_{2 x}^2 + f_{2 y}^2
\big]^{1/2}
}
{2}
\hspace{3cm}
\\[15pt]
&\hspace{3cm}
-
\bigg[
1 
+
\bigg(
\frac{f_{1 x} + f_{2 x}}{2}
\bigg)^2 
+ 
\bigg(
\frac{f_{1 y} + f_{2 y}}{2}
\bigg)^2
\hsx 
\bigg]^{1/2}
\bigg\}
\ \td x \hsy \td y
\ \geq \  
0
\end{align*}
\\[-.5cm]

or still, that
\\[-.5cm]

\begin{align*}
\iint\limits_Q \ 
\bigg\{
\bigg[
\bigg(
\frac{1}{2}
\bigg)^2
&+
\bigg(
\frac{f_{1 x}}{2}
\bigg)^2 
+ 
\bigg(
\frac{ f_{1 y}}{2}
\bigg)^2
\hsx 
\bigg]^{1/2}
\hsx + \hsx 
\bigg[
\bigg(
\frac{1}{2}
\bigg)^2
+
\bigg(
\frac{f_{2 x}}{2}
\bigg)^2 
+ 
\bigg(
\frac{ f_{2 y}}{2}
\bigg)^2
\hsx
\bigg]^{1/2}
\\[26pt] 
&
\hspace{1cm}
-
\bigg[
\bigg(
\frac{1}{2} + \frac{1}{2}
\bigg)
+
\bigg(
\frac{f_{1 x}}{2} + \frac{f_{2 x}}{2}
\bigg)^2
+
\bigg(
\frac{f_{1 y}}{2} + \frac{f_{2 y}}{2}
\bigg)^2
\hsx
\bigg]^{1/2}
\bigg\}
\ \td x \hsy \td y
\\[26pt]
&
\ \geq \  
0.
\end{align*}
\\[-.5cm]
\end{x}

\begin{x}{\small\bf LEMMA} \ 
\[
\sum\limits_{i = 1}^k \ 
(a_i^2 + b_i^2 + c_i^2)^{1/2} 
\ \geq \ 
\bigg[
\bigg(
\sum\limits_{i = 1}^k \ 
a_i
\bigg)^2
\hsx + \hsx 
\bigg(
\sum\limits_{i = 1}^k \ 
b_i
\bigg)^2
\hsx + \hsx 
\bigg(
\sum\limits_{i = 1}^k \ 
c_i
\bigg)^2
\hsx
\bigg]^{1/2}.
\]
\\[-.75cm]

To conclude that the foregoing integrand is nonnegative, take $k = 2$ and 
\[
\begin{cases}
\ 
\ds
a_1 = \frac{1}{2}, \quad b_1 = \frac{f_{1 x}}{2}, \quad c_1 = \frac{f_{1 y}}{2}
\\[15pt]
\ 
\ds
a_2 = \frac{1}{2}, \quad b_2 = \frac{f_{2 x}}{2}, \quad c_2 = \frac{f_{2 y}}{2}
\end{cases}
.
\]

Suppose that $f_1$ and $f_2$ are absolutely continuous in the sense of Tonelli and that 
equality obtains in Steiner $-$then 
the claim is that $f_1 - f_2$ is a constant.  
To establish this, observe first that

\[
\iint\limits_Q \ 
\{\ldots\} \ 
\td x \hsy \td y 
\ = \ 
0
\]
and since the integrand is nonnegative, it must be equal to zero almost everywhere in $Q$.  
This implies that

\[
f_{1 x} \ = \ f_{2 x}, 
\quad 
f_{1 y} \ = \ f_{2 y}
\]
almost everywhere in $Q$ or still, that 

\[
[(f_{1 x}  - f_{2 x})^2 + (f_{1 y}  - f_{2 y})^2 ]^{1/2} 
\ = \ 
0
\]
almost everywhere in $Q$.
\\[-.25cm]
\end{x}

\begin{x}{\small\bf NOTATION} \ 
$E \subset Q$ is the set consisting of 
\\[-.25cm]

\qquad 
(1) \ 
All lines $x = x_0$ such that $f_1(x_0, y)$, $f_2(x_0, y)$ are not both absolutely continuous in $y$.  
\\[-.5cm]

\qquad 
(2) \ 
All lines $y = y_0$ such that $f_1(x, y_0)$, $f_2(x, y_0)$ are not both absolutely continuous in $x$.  
\\[-.5cm]

\qquad 
(3) \ 
All points $(x, y)$ such that 

\[
f_{1 x} (x, y), 
\quad 
f_{2 x} (x, y), 
\quad
f_{1 y} (x, y), 
\quad 
f_{2 y} (x, y)
\]
are not all defined.
\\[-.5cm]

\qquad 
(4) \ 
All points $(x, y)$ at which

\[
[(f_{1 x}  - f_{2 x})^2 + (f_{1 y}  - f_{2 y})^2 ]^{1/2} 
\ \neq \ 
0.
\]
\\[-1cm]
\end{x}

\begin{x}{\small\bf \un{N.B.}}  \ 
$E$ has planar measure zero, hence for almost all points $(x_0, y_0) \in Q$ 
the lines $x = x_0$ and $y = y_0$ have in common with $E$ at most a set of linear measure zero. 
\\[-.25cm]
\end{x}

Fix one such point $(x_0, y_0)$ and let $(x, y)$ be any other point with the same property $-$then 

\[
\begin{cases}
\ 
f_1 (x, y) - f_1 (x_0, y_0) 
\ = \ 
\ds \int\limits_{x_0}^x \ 
f_{1 x} (x, y_0) \ \td x 
\hsx + \hsx
\int\limits_{y_0}^y \ 
f_{1 y} (x, y) \ \td y
\\[26pt]
\ 
f_2 (x, y) - f_2 (x_0, y_0) 
\ = \ 
\ds \int\limits_{x_0}^x \ 
f_{2 x} (x, y_0) \ \td x 
\hsx + \hsx
\int\limits_{y_0}^y \ 
f_{2 y} (x, y) \ \td y
\end{cases}
.
\]
Since apart from a set of linear measure zero the integrands on the right are equal, it thus follows that

\[
f_1 (x, y) - f_2 (x, y) 
\ = \ 
f_1 (x_0, y_0) - f_2 (x_0, y_0) , 
\]
which is true for almost all $(x, y)$ in $Q$, hence for all $(x, y)$ in $Q$ ($f_1$ and $f_2$ being continuous).
\\[-.5cm]

\begin{x}{\small\bf EXAMPLE} \ 
It can happen that equality prevails in Steiner, yet neither $f_1$ nor $f_2$ is ACT.
\\[-.5cm]

[Let $\phi(x)$ be a continuous monotonically increasing function such that $\phi^\prime (x) = 0$ almost everywhere 
and $\phi (0) = 0$, $\phi (1) = 1$.  
Working in $[0,2] \hsx \times \hsx [0,2]$, put

\[
\begin{cases}
\  f_1 (x, y) = 0 \hspace{1.8cm} (0 \leq x \leq 1, 0 \leq y \leq 2)
\\[4pt]
\  f_1 (x, y) = \phi(x - 1) \hspace{0.5cm} (1 \leq x \leq 2, 0 \leq y \leq 2)
\end{cases}
\]
and

\[
\begin{cases}
\  f_2 (x, y) = \phi (x) \hspace{1.1cm} (0 \leq x \leq 1, 0 \leq y \leq 2)
\\[4pt]
\  f_2 (x, y) = 1 \hspace{1.7cm} (1 \leq x \leq 2, 0 \leq y \leq 2)
\end{cases}
.
\]
Then

\[
L_Q [(f_1 + f_2)/2]  \ = \ 6, 
\quad
\begin{cases}
L_Q [f_1] \ = \ 6
\\[4pt]
L_Q [f_2] \ = \ 6
\end{cases}
\]

\qquad 
$\implies$

\[
6 
\ = \ 
\frac{6 + 6}{2} 
\ = \ \frac{12}{2}
\ = \ 
6.]
\]
\end{x}

\chapter{
$\boldsymbol{\S}$\textbf{9}.\quad EXTENSION PRINCIPLES}
\setlength\parindent{2em}
\setcounter{theoremn}{0}
\renewcommand{\thepage}{\S9-\arabic{page}}

\qquad 
Let $\phi : \sR \ra \R_{\geq 0}$ be a nonnegative rectangle function.
\\

\begin{x}{\small\bf PROBLEM} \ 
Determine conditions on $\phi$ which imply that $\phi$ can be extended to a measure on $\sB(Q)$ 
(the $\sigma$-algebra of Borel subsets of $Q$).
\\[-.5cm]
\end{x}

\begin{x}{\small\bf DEFINITION} \ 
$\phi$ satisfies \un{condition C} if for every choice of the systems

\[
\begin{cases}
\ 
r_1, \ldots, r_k
\\[4pt]
\ 
R_1, \ldots, R_n, \ldots
\end{cases}
\]
of oriented rectangles such that 

\[
r_i \cap r_j = \emptyset 
\quad (i \neq j)
\]
and 

\[
r_1 \cup \cdots \cup r_k 
\subset 
R_1 \cup \cdots \cup R_n \cup \cdots 
\quad 
\text{(finite or infinite)}
\]
there follows 

\[
\phi (r_1) + \cdots + \phi (r_k) 
\ \leq \ 
\phi (R_1) + \cdots + \phi (R_n) + \cdots \hsx .
\]
\\[-1.25cm]
\end{x}

\begin{x}{\small\bf DEFINITION} \ 
$\phi$ is \un{continuous} if for every $\varepsilon > 0$ there exists $\delta > 0$ such that 
$\phi(R) < \varepsilon$ for every oriented rectangle $R$ such that $\abs{R} < \delta$.
\\[-.5cm]
\end{x}

\begin{x}{\small\bf CRITERION} \ 
If $\phi$ is finitely additive and continuous, then $\phi$ satisfies condition C.
\\[-.5cm]
\end{x}

\begin{x}{\small\bf \un{N.B.}} \ 
Suppose that $\Phi$ is a Borel measure $-$then the restriction $\phi \equiv \restr{\Phi}{\sR}$ satisfies condition C.
\\[-.5cm]


[Put 
\\[-1.25cm]
\begin{align*}
&
B_1 = R_1,
\\[4pt]
&
B_2 = R_2 \backslash R_1,
\\[4pt]
&
\hspace{0.5cm}
\vdots
\\[4pt]
&
B_n = R_n \backslash (R_1 \cup \cdots \cup R_{n - 1}),
\\[4pt]
&
\hspace{0.5cm}
\vdots
\ .
\end{align*}
Then
\begin{align*}
\phi (r_1) + \cdots + \phi (r_k) \ 
&=\ 
\Phi (r_1) + \cdots + \Phi (r_k) 
\\[11pt]
&=\ 
\Phi (r_1 \cup \cdots \cup r_k)
\\[11pt]
&\leq\ 
\Phi (R_1 \cup \cdots \cup R_n \cup \cdots)
\\[11pt]
&=\ 
\Phi (B_1 \cup \cdots \cup B_n \cup \cdots)
\\[11pt]
&=\ 
\Phi (B_1) + \cdots + \Phi (B_n) + \cdots
\\[11pt]
&\leq\ 
\Phi (R_1) + \cdots + \Phi (R_n) + \cdots
\\[11pt]
&=\ 
\phi (R_1) + \cdots + \phi (R_n) + \cdots \hsx .]
\end{align*}
\\[-1.25cm]
\end{x}

\begin{x}{\small\bf NOTATION} \ 
Given a set $E \subset Q$, let
\[
\Gamma^* (\phi ; E)
\ = \ 
\inf \ \sum \ \phi(R_n), 
\]
where the inf is taken over all rectangles $R_1, \ldots, R_n, \ldots$ (finite or infinite) of oriented 
rectangles in $Q$ such that $E \subset \bigcup R_n$ 
(take $\Gamma(\emptyset, \phi) = 0$).
\\[-.5cm]
\end{x}

\begin{x}{\small\bf LEMMA} \ 
Suppose that $\phi$ satisfies condition C 
$-$then  $\Gamma^* (\phi ; -)$ is a metric outer measure.
\\[-.5cm]
\end{x}

\begin{x}{\small\bf NOTATION} \ 
Put 
\[
\Gamma (\phi ; -) 
\ = \ 
\restr{\Gamma^* (\phi ; -)}{\sB(Q)},
\]
a measure on $\sB (Q)$.
\\[-.5cm]
\end{x}

\begin{x}{\small\bf THEOREM} \ 
$\phi$ extends to a measure on $\sB (Q)$ iff $\phi$ satisfies condition C.
\\[-.25cm]

PROOF \ 
The necessity follows from \#5 and the sufficiency follows from \#7 
(obviously, $\forall \ R \in \sR$, $\Gamma (\phi ; R) = \phi(R)$).
\\[-.5cm]
\end{x}

\begin{x}{\small\bf LEMMA} \ 
If $\Phi$ and $\Psi$ are Borel measures and if $\Phi(R) = \Psi (R)$ $(\forall \ R \in \sR)$, 
then $\Phi(E) = \Psi (E)$ $(\forall \ E \in \sB (Q)).$
\\

Suppose that $f: Q \ra \R$ is of bounded variation in the sense of Tonelli and recall that 

\[
\begin{cases}
\ 
W_x (f; R) 
\ = \ 
\ds
\int\limits_c^d \ 
\tV_x (f; y) 
\ \td y
\\[26pt]
\ 
W_y (f; R) 
\ = \ 
\ds
\int\limits_a^b \ 
\tV_y (f; x) 
\ \td x
\end{cases}
.
\]

It is clear that

\[
\begin{cases}
\ 
W_x (f; -) 
\\[8pt]
\ 
W_y (f; -) 
\end{cases}
\]
are finitely additive and it can be shown that they are continuous.  
Therefore

\[
\begin{cases}
\ 
W_x (f; -) 
\\[8pt]
\ 
W_y (f; -) 
\end{cases}
\]
satisfy condition C (cf. \# 4), thus they each admit a unique extension to a measure on $\sB(Q)$, denoted 

\[
E \ra \ 
\begin{cases}
\ 
W_x (f; E) 
\\[8pt]
\ 
W_y (f; E) 
\end{cases}
\quad (E \in \sB(Q)).
\]
Accordingly there are Lebesgue decompositions

\[
\begin{cases}
\ 
\ds
W_x (f; E) 
\ = \ 
\iint\limits_E \ 
\abs{f_x} \ \td \Lm^2 
\hsx + \hsx 
W_x^0 (f; E)
\\[26pt]
\ 
\ds
W_y (f; E) 
\ = \ 
\iint\limits_E \ 
\abs{f_y} \ \td \Lm^2 
\hsx + \hsx 
W_y^0 (f; E)
\end{cases}
,
\]
where

\[
\begin{cases}
\ 
W_x^0 (f; -) 
\\[8pt]
\ 
W_y^0 (f; -) 
\end{cases}
\]
are singular.
\end{x}

\chapter{
$\boldsymbol{\S}$\textbf{10}.\quad ONE VARIABLE REVIEW}
\setlength\parindent{2em}
\setcounter{theoremn}{0}
\renewcommand{\thepage}{\S10-\arabic{page}}

\qquad
In the Fr\'echet process, take for $\X$ the quasi linear functions $\Gamma$ on $[0,1]$, 
take for $d$ the metric defined by the prescription

\[
d(\Gamma_1,\Gamma_2) 
\ = \ 
\int\limits_0^1 \abs{\Gamma_1(x) - \Gamma_2(x)} \ \tdx,
\]
and take for $F$ the elementary length $-$then the lower semicontinuity is manifest, as is property (A).  
Here $\ov{\X} = \Lp^1[0,1]$ and property (B) is satisfied.
\\[-.25cm]

\begin{x}{\small\bf DEFINITION} \ 
Put
\[
\text{\Large\textnormal{\textturnh}}[f] 
\ = \ 
\bar{F}(f) \qquad (f \in \Lp^1[0,1])
\]
and call it the \un{generalized variation} of $f$.
\\[-.25cm]
\end{x}

\begin{x}{\small\bf DEFINITION} \ (gBV) \ 
A function $f \in \Lp^1[0,1]$ is of \un{generalized bounded} \un{variation} if
\[
\text{\Large\textnormal{\textturnh}}[f] \ < \ +\infty.
\]
\\[-1cm]
\end{x}

\begin{x}{\small\bf NOTATION} \ 
$\gBV[0,1]$ is the set of all functions of generalized bounded variation.
\\[-.5cm]
\end{x}

\begin{x}{\small\bf THEOREM} \ 
Let $f \in \Lp^1[0,1]$ $-$then $f$ is of generalized bounded variation iff there is a 
$g \in\ \Lp^1[0,1]$ which is equal almost everywhere to $f$ and $T_g[0,1] < +\infty$.
\\[-.5cm]

Therefore
\[
\BV[0,1] \ \subset \ \gBV[0,1].
\]
\\[-1.25cm]
\end{x}

\begin{x}{\small\bf THEOREM} \ 
Suppose that $f \in \gBV[0,1]$ $-$then
\[
\text{\Large\textnormal{\textturnh}}[f] \ = \ \inf \{T_g[0,1]: g = f \text{ almost everywhere}\}.
\]
\\[-1.5cm]
\end{x}


\begin{x}{\small\bf RAPPEL} \ 
Given an $f \in \Lp^1[0,1]$, $C_\ap (f)$ is the set of points of approximate continuity.
\\
\end{x}

\begin{x}{\small\bf \un{N.B.}}  \ 
$C_\ap(f)$ is a subset of $[0,1]$ of full measure.
\\[-.25cm]
\end{x}

\begin{x}{\small\bf LEMMA} \ 
If $f \in \Lp^1[0,1]$, then
\[
\text{\Large\textnormal{\textturnh}}[f] \ = \ \sup \ \sum\limits_{i=1}^{n-1} \ \abs{f(x_{i+1}) - f(x_i)}, 
\]
where the supremum is taken over all finite collections of points $x_i \in C_\ap(f)$ subject to $x_i < x_{i+1}$.
\\[-.5cm]

[Note: \ 
If $E \subset C_\ap(f)$ is a subset of full measure, then the supremum can be taken over the $x_i \in E$.]
\\[-.25cm]
\end{x}

\begin{x}{\small\bf RAPPEL} \ 
If $f_n \ra f$ in $\Lp^1[0,1]$, then there is a subsequence $\{f_{n_k}\}$ such that $\{f_{n_k}\} \ra f$ almost everywhere.
\\[-.25cm]
\end{x}

\begin{x}{\small\bf LEMMA} \ 
{\textnormal{\Large\textturnh}} is lower semicontinuous w.r.t. convergence almost everywhere, 
i.e., if $f_1, f_2, \ldots$ is a sequence in $\Lp^1[0,1]$ that converges almost everywhere to $f \in \Lp^1[0,1]$, then
\[
\text{\Large\textnormal{\textturnh}}[f] \ 
\leq \ 
\liminf\limits_{n \ra \infty} \  \text{\Large{\textnormal{\textturnh}}}[f_n].
\]
\\[-1cm]
\end{x}

\begin{x}{\small\bf DEFINITION} \ 
The \un{essential derivative} of $f$ at a point $x$ is the derivative of $f$ 
computed at $x$ after deleting a set of Lebesgue measure 0.
\\[-.25cm]
\end{x}

\begin{x}{\small\bf THEOREM} \ 
Suppose that $\text{\Large\textnormal{\textturnh}}[f]$ is finite 
$-$then the \un{essential derivative} of $f$, denoted still by $f^\prime$, exists almost everywhere and
\[
\text{\Large\textnormal{\textturnh}}[f] 
\ \geq \ 
\int\limits_0^1 \abs{f^\prime} dx.
\]
Moreover equality obtains iff $f$ is equivalent to an absolutely continuous function. 
\end{x}

\chapter{
$\boldsymbol{\S}$\textbf{11}.\quad EXTENDED LEBESGUE AREA}
\setlength\parindent{2em}
\setcounter{theoremn}{0}
\renewcommand{\thepage}{\S11-\arabic{page}}

\qquad
In the Fr\'echet process, take for $\X$ the quasi linear functions $\Pi$ on $[0,1] \times [0,1]$ $(= Q)$, 
take for $d$ the metric defined by the prescription
\[
d(\Pi_1,\Pi_2) \ = \ \iint\limits_Q \abs{\Pi_1(x, y) - \Pi_2(x, y)} \ \tdx \hsy \tdy,
\]
and take for $F$ the elementary area $-$then lower semicontinuity is manifest, as is property (A).  
Here $\ov{\X} = \Lp^1(Q)$ and property (B) is satisfied.
\\[-.25cm]

\begin{x}{\small\bf DEFINITION} \ 
Put
\[
\text{\Large\textnormal{\textturnh}}_Q[f] \ = \ \ovs{F}(f) \qquad (f \in \Lp^1 (Q))
\]
and call it the \underline{generalized variation} of $f$.
\\[-.25cm]
\end{x}

\begin{x}{\small\bf EXTENSION PRINCIPLE} \ 
Suppose that $f:Q \ra \R$ is continuous $-$then
\[
\text{\Large\textnormal{\textturnh}}_Q[f] \ = \ L_Q[f].
\]

\end{x}

\begin{x}{\small\bf \un{N.B.}} \ 
Therefore $\text{\Large\textnormal{\textturnh}}_Q$ can be viewed as an ``area functional'' on $\Lp^1(Q)$, 
there being no a priori assumption of continuity, which justifies calling 
$\text{\Large\textnormal{\textturnh}}_Q$ 
extended Lebesgue area.
\end{x}

\begin{x}{\small\bf LEMMA} \ 
Suppose that $f : Q \ra \R$ is continuous.
\\[-.25cm]

\qquad
\textbullet \quad
If $\Lm_Q [f] < +\infty$, then for every $\varepsilon > 0$ there is a $\delta > 0$ such that if 
$g : Q \ra \R$ is continuous and 
\[
\abs{f(x, y) - g(x, y)} 
\ < \ 
\delta
\]
on a set of measure greater than $1 - \delta$, then 
\[
\Lm_Q [g]
\ > \ 
\Lm_Q [f] - \varepsilon.
\]

\qquad
\textbullet \quad
If $\Lm_Q [f] = +\infty$, then for every $M > 0$ there is a $\delta > 0$ such that if 
$g : Q \ra \R$ is continuous and 
\[
\abs{f(x, y) - g(x, y)} 
\ < \ 
\delta
\]
on a set of measure greater than $1 - \delta$, then 
\[
\Lm_Q [g]
\ > \ 
M.
\]
\\[-1.75cm]
\end{x}

There are two possibilities: 
\[
\Lm_Q [f] 
\ < \ 
+\infty
\quad \text{or} \quad 
\Lm_Q [f] 
\ =\ 
+\infty.
\]
For sake of argument, consider the first of these.
\\[-.5cm]

Since uniform convergence of $\{\Pi_n( x, y)\}$ to $f(x, y)$ implies that 
$\td (\Pi_n, f)$ converges to zero, 
it follows that 
$\text{\Large\textnormal{\textturnh}}_Q [f] \leq \Lm_Q [f]$.  
To go the other way, take $\varepsilon > 0$, let  $\delta > 0$ be per supra, 
and choose a quasi linear function $\Pi$ such that 
\[
\iint\limits_Q \ \abs{f - \Pi} \ \td \Lm^2 
\ < \ 
\delta^2.
\]
Then 
\[
\abs{f(x, y) - \Pi (x, y)} 
\ < \ 
\delta
\]
on a set of measure greater than $1 - \delta$, hence
\[
\Lm_Q [\Pi]
\ > \ 
\Lm_Q [f] - \varepsilon
\]
\\[-1.75cm]

\qquad $\implies$
\[
\text{\Large\textnormal{\textturnh}}_Q[f]
\ \geq \ 
\Lm_Q [f] - \varepsilon.
\]
\\[-1cm]

There is also a Ge\"ocze version of these considerations.
\\[-.5cm]

\begin{x}{\small\bf DEFINITION} \ 
Let $f \in \Lp^1 (Q)$ and let $R \subset Q$ be an oriented rectangle, 
thus in the usual notation, 

\vspace{-.25cm}
\[
\begin{cases}
\ a \leq x \leq b \hspace{0.5cm} (a < b)
\\[8pt]
\ c \leq y \leq d \hspace{0.5cm} (c < d)
\end{cases}
, \ 
\abs{R} \hsx = \hsx (b - a) \hsy (d - c).
\]
Then $R$ is said to be \un{admissible} if 
$f(x, y)$ is approximately continuous in $x$ for almost all $y$ 
on the boundary lines of $R$ parallel to the $y$ axis and if 
$f(x, y)$ is approximately continuous in $y$ for almost all $x$ 
on the boundary lines of $R$ parallel to the $x$ axis.
\\[-.5cm]

[Note: \ 
A subdivision $D$ of $Q$ into nonoverlapping oriented rectangles 
$R$ is admissible provided this is the case of each of the $R$.] 
\\[-.5cm]
\end{x}

Using this data, one can arrive at the extended Ge\"ocze area, denoted

\[
\extgeocze_Q [f].
\]
\\[-1.25cm]

\begin{x}{\small\bf THEOREM} \ 
\[
\extgeocze_Q (f)
\ = \ 
\text{\Large\textnormal{\textturnh}}_Q[f].
\]
\\[-1cm]
\end{x}

\begin{x}{\small\bf \un{N.B.}} \ 
Recall that 

\[
\Gamma_Q [f] 
\ = \ 
\Lm_Q [f] 
\quad 
(f \in C(Q)), 
\]
i.e., 
\[
\text{Ge\"ocze area = Lebesgue area.}
\]
\end{x}

\chapter{
$\boldsymbol{\S}$\textbf{12}.\quad THEORETICAL SUMMARY}
\setlength\parindent{2em}
\setcounter{theoremn}{0}
\renewcommand{\thepage}{\S12-\arabic{page}}

\qquad 
What is said below for the integrable case runs parallel to what has been said for the continuous case. 
\\[-.25cm]

\begin{x}{\small\bf DEFINITION} \ (gBVT) \ 
Let $f \in \Lp^1 (Q)$ $-$then $f$ is said to be of 
\un{generalized bounded variation in the sense of Tonelli} if

\[
\begin{cases}
\ \ds 
\int_0^1 \ \text{\Large\textnormal{\textturnh}}[f(-,y)] \ \tdy \ < \ +\infty 
\\[26pt]
\ \ds 
\int_0^1 \ \text{\Large\textnormal{\textturnh}}[f(x,-)] \ \tdx \ < \ +\infty
\end{cases}
.
\]
\end{x}

The gBVT-functions can be characterized.
\\[-.25cm]

\begin{x}{\small\bf THEOREM} \ 
Let $f \in \Lp^1 (Q)$ $-$then $f$ is of generalized bounded variation in the sense of Tonelli 
iff there are functions $g$ and $h$ equal to $f$ almost everywhere in $Q$ such that
\[
\begin{cases}
\ \ds 
\ \int_0^1 \ \textnormal{V}_x (g;y) \ \tdy \ < \ +\infty 
\\[26pt]
\ \ds 
\int_0^1 \ \textnormal{V}_y (g;x) \ \tdx \ < \ +\infty
\end{cases}
.
\]
\\[-1cm]
\end{x}

\begin{x}{\small\bf REMARK} \ 
Suppose that $f$ is gBVT $-$then it can be shown that there is a function $k$ equal to $f$ almost everywhere in $Q$ such that

\[
\begin{cases}
\ \ds 
\int_0^1 \ \tV_x (k;y) \ \tdy \ < \ +\infty 
\\[26pt]
\ \ds 
\ \int_0^1 \ \tV_y (k;x) \ \tdx \ < \ +\infty
\end{cases}
.
\]
\end{x}

\begin{x}{\small\bf \un{N.B.}} \ 
\[
f \ \BVT \ \implies \  f \ \gBVT.
\]

[Note: \
 Recall that $f$ \BVT means, in particular, that $f \in C(Q)$, hence $f \in \Lp^1(Q).]$
 \\[-.25cm]
\end{x}

\begin{x}{\small\bf THEOREM} \ 
$\text{\Large\textnormal{\textturnh}}_Q[f] \ < \ +\infty$ iff $f$ is gBVT.
\\[-.25cm]
\end{x}

\begin{x}{\small\bf THEOREM} \ 
Suppose that $f$ is gBVT $-$then the essential partial derivatives $f_x$ and $f_y$ exist almost everywhere, are integrable, and
\[
\text{\Large\textnormal{\textturnh}}_Q[f] 
\ \geq \ 
\iint\limits_Q \ 
\left[1 + f_x^2 + f_y^2\right]^{1/2} 
\ \td x \hsy \td y.
\]
\\[-1cm]
\end{x}

\begin{x}{\small\bf DEFINITION} \ (gACT) \ 
Suppose that $f$ is gBVT $-$then $f$ is said to be 
\un{generalized absolutely continuous in the sense of Tonelli} if $f$ coincides almost everywhere 
with a function $g$ which is absolutely continuous w.r.t. $x$ for almost all $y$ 
and absolutely continuous w.r.t $y$  for almost all $x$.
\\[-.25cm]
\end{x}

\begin{x}{\small\bf SCHOLIUM} \ 
\\[-.25cm]

\qquad
\textbullet \ 
If $f$ is gBVT and if 

\[
\text{\Large\textnormal{\textturnh}}_Q[f]
\ = \ 
\iint\limits_Q \ 
\left[1 + f_x^2 + f_y^2\right]^{1/2} 
\ \td x \hsy \td y,
\]
then $f$ is gACT.
\\[-.25cm]

\qquad
\textbullet \ 
If $f$ is gACT, then

\[
\text{\Large\textnormal{\textturnh}}_Q[f]
\ = \ 
\iint\limits_Q \ 
\left[1 + f_x^2 + f_y^2\right]^{1/2} 
\ \td x \hsy \td y.
\]
\end{x}

\chapter{
$\boldsymbol{\S}$\textbf{13}.\quad VARIANTS}
\setlength\parindent{2em}
\setcounter{theoremn}{0}
\renewcommand{\thepage}{\S13-\arabic{page}}

\qquad 
Up to this point, the discussion has taken

\[
Q 
\ = \ 
[0,1] \hsx \times \hsx [0,1]
\]
as the domain of discourse.  
Of course, matters can be extended with little change when $Q$ is replaced by 
\[
[a,b] \hsx \times \hsx [c,d].
\]
This done, the next step is to replace $Q$ by a nonempty open subset $\Omega \subset \R^2$.
\\[-.25cm]

\begin{x}{\small\bf RAPPEL} \ 
A continuous function $f : \hsx ]a,b[ \hsx  \ra \R$ is of bounded variation in a nonempty open interval $]a,b[ \hsx \subset \R$ provided
\[
T_f \hsx ]a,b[ 
\ < \ 
+\infty.
\]
\\[-1.5cm]
\end{x}

\begin{x}{\small\bf DEFINITION} \ 
A continuous function $f : \Omega \ra \R$ is of bounded variation in a nonempty open subset $\Omega \subset \R$ provided

\[
T_f \Omega 
\ < \ 
+\infty,
\]
where
\[
T_f \Omega 
\ = \ 
\sum\limits_n \ 
T_f \hsx ]a_n, b_n[\hsx, 
\]
the nonempty open intervals $]a_n, b_n[$ running through the connected components of $\Omega$ (admit $\pm\infty$).
\\[-.25cm]
\end{x}

\begin{x}{\small\bf NOTATION} \ 
Let $\Omega$ be a nonempty open subset of $\R^2$.
\\[-.5cm]

\qquad 
\textbullet \ 
For any real number $\bar{x}$, let $\Omega(\bar{x})$ denote the open linear set which is 
the intersection of $\Omega$ with the straight line $x = \bar{x}$. 
\\[-.5cm]

\qquad 
\textbullet \ 
For any real number $\bar{y}$, let $\Omega(\bar{y})$ denote the open linear set which is 
the intersection of $\Omega$ with the straight line $y = \bar{y}$. 
\\[-.25cm]
\end{x}

Given a continuous function $f : \Omega \ra \R$, introduce 

\[
\begin{cases}
\ 
\tV_x (f; \bar{y}; \Omega) 
\ = \ 
\tT_f \Omega (\bar{y})
\\[8pt]
\
\tV_y (f; \bar{x}; \Omega)
\ = \ 
\tT_f \Omega (\bar{x})
\end{cases}
.
\]

[Note: \ 
Take 

\[
\begin{cases}
\ 
\tV_x = 0 \quad \text{if} \quad \Omega(\bar{y}) = \emptyset
\\[8pt]
\ 
\tV_y = 0 \quad \text{if} \quad \Omega(\bar{x}) = \emptyset
\end{cases}
.]
\]
\\[-1cm]

\begin{x}{\small\bf LEMMA}\ 

\[
\begin{cases}
\ 
\tV_x (f; \bar{y}; \Omega)
\ 
\text{is a lower semicontinuous function of $\bar{y}$}
\\[8pt]
\ 
\tV_y (f; \bar{x}; \Omega)
\ 
\text{is a lower semicontinuous function of $\bar{x}$}
\end{cases}
\
\text{in $]-\infty, +\infty[$}.
\]
\\[-1cm]
\end{x}

\begin{x}{\small\bf DEFINITION} \ (BVT) \ 
$f$ is said to be of \un{bounded variation in the} \un{sense of Tonelli} if

\[
\begin{cases}
\ \ds
\int\limits_{-\infty}^{+\infty} \ 
\tV_x (f; \bar{y}; \Omega) \ 
\td \bar{y} 
\ < \ 
+\infty
\\[26pt]
\ \ds
\int\limits_{-\infty}^{+\infty} \ 
\tV_y (f; \bar{x}; \Omega) \ 
\td \bar{x} 
\ < \ 
+\infty
\end{cases}
.
\]
\\[-1cm]
\end{x}

\begin{x}{\small\bf LEMMA} \ 
Suppose that $f : \Omega \ra \R$ is of bounded variation in the sense of Tonelli $-$then

\[
\begin{cases}
\ \ds
f_x \ = \ \frac{\partial f}{\partial x}
\\[15pt]
\ \ds
f_y \ = \ \frac{\partial f}{\partial y}
\end{cases}
\ \text{exists almost everywhere in $\Omega$}
\]
and

\[
\begin{cases}
\ \ds
\iint\limits_\Omega \ 
\abs{f_x (x, y)} \ 
\tdx \hsy \tdy
\ \leq \ 
\int\limits_{-\infty}^{+\infty} \ 
\tV_x (f; \bar{y}; \Omega) \ 
\td \bar{y} 
\ < \ 
+\infty
\\[26pt]
\ \ds
\iint\limits_\Omega \ 
\abs{f_y (x, y)} \ 
\tdx \hsy \tdy
\ \leq \ 
\int\limits_{-\infty}^{+\infty} \ 
\tV_y (f; \bar{x}; \Omega) \ 
\td \bar{x} 
\ < \ 
+\infty
\end{cases}
\]
\\[-1cm]

\qquad 
$\implies$
\[
\begin{cases}
\ 
f_x 
\\[4pt]
f_y 
\end{cases}
\ \in \ \Lm^1 (\Omega).
\]
\\[-1.25cm]
\end{x}

Another setting for the theory is a nonempty open subset $\Omega \subset \R^2$, 
$\Lm^1 (Q)$ then being replaced by $\Lm^1 (\Omega)$, 
the analog of a gBVT function now being an element of BV$\Lm^1 \Omega$.
\\[-.5cm]

\begin{x}{\small\bf DEFINITION} \ 
Let $f \in \Lm^1 (\Omega)$ $-$then $f$ is a \un{function of bounded variation} in $\Omega$ 
if the distributional partial derivatives of $f$ are finite signed Radon measures
\\[-.5cm]

\[
\begin{cases}
\ \ds
\mu_x \hsx : \hsx  \int\limits_\Omega \ 
f \frac{\partial \phi}{\partial x} \ \td x
\ = \ 
-\int\limits_\Omega \ 
\phi \ \td \mu_x
\\[26pt]
\ \ds
\mu_y \hsx : \hsx  \int\limits_\Omega \ 
f \frac{\partial \phi}{\partial y} \ \td y
\ = \ 
-\int\limits_\Omega \ 
\phi \ \td \mu_y
\end{cases}
\ \forall \ \phi \in C_c^\infty (\Omega)
\]
of finite total variation. 
\\[-.5cm]
\end{x}

\begin{x}{\small\bf NOTATION} \ 
$\BV\Lm^1\Omega$ is the set of functions of bounded variation in $\Omega$.  
\\[-.5cm]
\end{x}

Given $g \in \Lm^1 (\Omega)$, put
\[
\tV_T (g; \Omega) 
\ = \ 
\int\limits_{-\infty}^{+\infty} \ 
\tV_x (g; \bar{y}; \Omega) 
\ \td \bar{y} 
\ + \hsx 
\int\limits_{-\infty}^{+\infty} \ 
\tV_y (g; \bar{x}; \Omega) 
\ \td \bar{x}.
\]
\\[-1cm]


\begin{x}{\small\bf THEOREM} \ 
Let $f \in \Lm^1 (\Omega)$ $-$then $f \in \BV\Lm^1 \Omega$ iff
\[
\inf \{\tV_T (g; \Omega) : g = f \ \text{almost everywhere} \} 
\ < \ 
+\infty.
\]
\end{x}

\newpage

\centerline{\textbf{\large REFERENCES}}
\setcounter{page}{1}
\setcounter{theoremn}{0}
\renewcommand{\thepage}{References-\arabic{page}}
\vspace{0.75cm}

\[
\text{BOOKS}
\]

\begin{rf}
Cesari, Lamberto, \textit{Surface Area}, Annals of Mathematics Studies, Number 35, Princeton University Press, 1956.
\end{rf}

\begin{rf}
Evans, L. C. and Gariepy, R. F., \textit{Measure theory and fine properties of functions}, CRC Press, Boca Raton, 1992.
\end{rf}

\begin{rf}
Federer, Herbert, \textit{Geometric Measure Theory}, Springer-Verlag, New York, 1969.
\end{rf}

\begin{rf}
Goffman, Casper, Nishiura, Toga, and Waterman, Daniel, \textit{Homeomorphisms in Analysis}, American Mathematical Society, 1997.
\end{rf}

\begin{rf}
Rado,Tibor, \textit{Length and Area}, American Mathematical Society Colloquium Publications, 30, 1948.
\end{rf}

\begin{rf}
Rado,Tibor and Reichelderfer, Paul V., \textit{Continuous Transformations in Analysis: With an Introduction to Algebraic Topology}, 
Die Grundlehren der Mathematischen Wissenschaften, vol 75. Springer, Berlin, Heidelberg 1955.
\end{rf}

\begin{rf}
Saks, S.,  \textit{Theory of the integral}, Warszawa, 1937.
\end{rf}

\begin{rf}
Tonelli, L.,  \textit{Fondimenti di calcolo delle variazioni}, vols. 1 and 2, Bologna, Zanichelli, 1922.
\end{rf}


\setcounter{theoremn}{0}

\[
\text{ARTICLES}
\]

\begin{rf}
Banach, S., 
Sur les lignes rectifiables et les surfaces dout l'aire est fini, 
\textit{Fund. Math.} \textbf{7} (1925), 225-236.
\end{rf}

\begin{rf}
Besicovitch, A. S.,  
On the definition and value of the area of a surface, 
\textit{Quart. J. Math. Oxford Ser.} \textbf{16} (1945), 86-102.
\end{rf}

\begin{rf}
Blumberg, H., 
New properties of all real functions, 
\textit{Trans. Amer. Math. Soc.} \textbf{24} (1922), 113-128.
\end{rf}

\begin{rf}
Breckenridge, J. C., and Nishiura, T., 
Two examples in surface area theory, \textit{Michigan Math. J.} \textbf{19} (1972), 157-160.
\end{rf}

\begin{rf}
Bruckner, A. M., and Goffman, C., 
Differentiability through change of variables,  \textit{Proc. Amer. Math. Soc.} \textbf{61} (1976), 235-241.
\end{rf}

\begin{rf}
Cesari, L., 
Area and Represenation of Surfaces, \textit{Bull. Amer. Math. Soc.} vol. \textbf{56} (1950), 218-232.
\end{rf}

\begin{rf}
Cesari, L., 
Recent results in surface area theory, \textit{Amer. Math. Monthly}  \textbf{66}, (1959), 173-192.
\end{rf}

\begin{rf}
Federer, H.,  
Surface Area (I), (II) \textit{Trans. Amer. Math. Soc.}  \textbf{55},  (1944), 420-456.
\end{rf}

\begin{rf}
Federer, H.,  
Measure and Area, \textit{Bull. Amer. Math. Soc.} \textbf{58} (1952), 306-378.
\end{rf}

\begin{rf}
Federer, H.,  
The area of a nonparametric surface, \textit{Proc. Amer. Math. Soc.} vol. \textbf{11} (1960), 436-439.
\end{rf}

\begin{rf}
Fleming, W. H., 
Functions whose partial derivatives are measures, \textit{Illinois J. Math.} \textbf{4} (1960), 452-478.
\end{rf}

\begin{rf}
Hughs, R. E., 
Functions of BVC type, \textit{Proc. Amer. Math. Soc.} \textbf{12} (1961), 698-701.
\end{rf}

\begin{rf}
Reifenberg, R. E., 
Parametric surfaces (I), (II), \textit{Proc. Cambridge Philos. Soc.} vol. \textbf{47}  (1951), 687-698.
\end{rf}

\begin{rf}
Tall F. D., 
The density topology, \textit{Pacific J. Math.}  \textbf{62}, (1976), 275-284.
\end{rf}

\begin{rf}
Ward, D. J., 
A counterexample in area theory, \textit{Proc. Cambridge Philos. Soc.} vol. \textbf{60} (1964), 821-845.
\end{rf}

\begin{rf}
Youngs, J. W. T., 
Curves and Surfaces, \textit{Amer. Math. Monthly} \textbf{51} (1944), 1-11.
\end{rf}

\begin{rf}
Zahorski, Z., 
Sur la premiere derivee, \textit{Trans. Amer. Math. Soc.} \textbf{69} (1950), 1-54.
\end{rf}